\documentclass{amsart}
\usepackage{amssymb}
\usepackage{amsfonts}
\usepackage{diagrams}
\usepackage[dvips]{graphicx}
\usepackage{float}

\newtheorem{theorem}{Theorem}[section]
\theoremstyle{plain}

\newtheorem{corollary}[theorem]{Corollary}

\newtheorem{definition}[theorem]{Definition}

\newtheorem{lemma}[theorem]{Lemma}

\newtheorem{proposition}[theorem]{Proposition}
\newtheorem{remark}[theorem]{Remark}

\newtheorem*{main}{Main Theorem}
\numberwithin{equation}{section} \DeclareMathOperator{\nbhd}{nbhd}
\DeclareMathOperator{\Hom}{Hom} \DeclareMathOperator{\im}{im}
\DeclareMathOperator{\coker}{coker}

\begin{document}
\title[Higher-Order Linking Forms for Knots]{Higher-Order Linking Forms for
Knots}
\author{Constance Leidy}
\email{leidy@math.rice.edu}

\begin{abstract}
We construct examples of knots that have isomorphic $n$th-order
Alexander modules, but non-isomorphic $n$th-order linking forms,
showing that the linking forms provide more information than the
modules alone. This generalizes work of Trotter \cite{trotter}, who
found examples of knots that have isomorphic classical Alexander
modules, but non-isomorphic classical Blanchfield linking forms.

\end{abstract}

\maketitle

\section{Introduction}

In 1973, Trotter \cite{trotter} found examples of knots that have
isomorphic classical Alexander modules, but non-isomorphic classical
Blanchfield linking forms. Recently, T. Cochran \cite{nckt}\ defined
higher-order Alexander modules, $\mathcal{A}_{n}\left( K\right) $,
of a knot, $K$, and higher-order linking forms, $\mathcal{B\ell
}_{n}\left( K\right) $, which are linking forms defined on
$\mathcal{A}_{n}\left( K\right) $. When $n=0$, these invariants are
just the classical Alexander module and Blanchfield linking form.
The question was posed in \cite{nckt} whether Trotter's result
generalized to the higher-order invariants. We show that it does.
The following is our main theorem.

\begin{main}
For each $n\geq 0$, there exist knots $K_{0}$ and $K_{1}$ such that
$ \mathcal{A}_{i}\left( K_{0}\right) \cong \mathcal{A}_{i}\left(
K_{1}\right) $ for $0\leq i\leq n$ and $\mathcal{B\ell }_{i}\left(
K_{0}\right) \cong \mathcal{B\ell }_{i}\left( K_{1}\right) $ for
$0\leq i<n$, but $\mathcal{ B\ell }_{n}\left( K_{0}\right) \ncong
\mathcal{B\ell }_{n}\left( K_{1}\right) $.
\end{main}

When $n=1$, a particular example of the main theorem is the
following pair of knots.  The construction of them will be explained
later in this paper.

\begin{figure}[H]
\centering
\includegraphics[scale=.75]{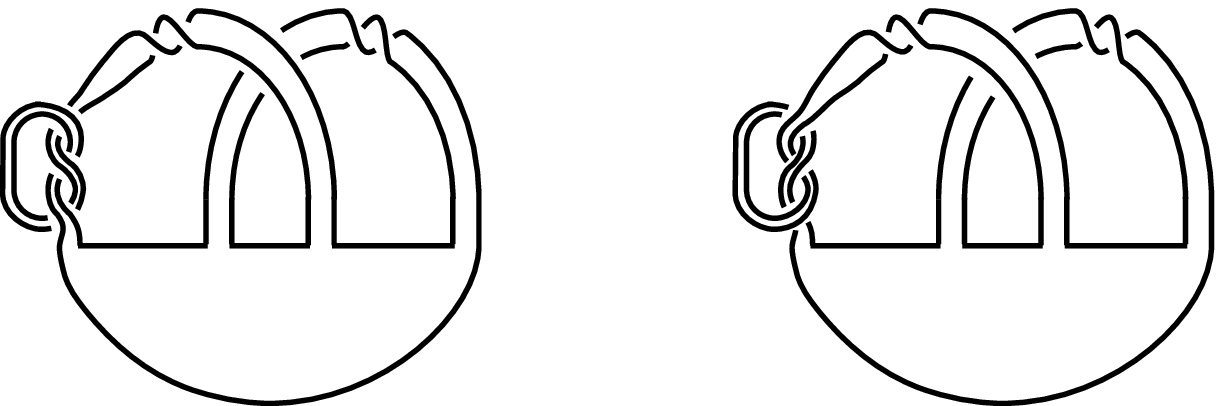}
\end{figure}

We shall work with classical, oriented knots in the PL category. We
now review some notions of classical knot theory. We refer the
reader to \cite{camerongordon}, \cite{lickorish}, and \cite{rolfsen}
as knot theory resources.  Recall that by Alexander duality, the
$p$-th reduced homology of the exterior of the knot is trivial
except when $p=1$, in which case it is $\mathbb{Z}$, generated by
the meridian. It follows that $\frac{G}{ G^{\prime }}\cong
\mathbb{Z}$, where $G$ is the fundamental group of the exterior.
Hence, we can take the infinite cyclic cover of the exterior.  The
classical Alexander module of a knot is defined to be the first
homology of this infinite cyclic cover of the exterior of the knot,
considered as a $ \mathbb{Z}\left[ t,t^{-1}\right] $-module. Here
the module structure results from the action $x\ast t=\mu
^{-1}x\mu$, where $\mu$ is the meridian of the knot. Furthermore,
since the fundamental group of the infinite cyclic cover is the
commutator subgroup, $G^{\prime }$, it follows that the Alexander
module is simply $\frac{ G^{\prime }}{G^{\prime \prime }}$
considered as a right $\mathbb{Z}\left[ \frac{G}{G^{\prime }}\right]
$-module.

A Seifert surface for a knot, $K$, is a connected, bicollared,
compact surface in $S^{3}$ whose boundary is $K$.  For a choice of
Seifert surface, $F$, and bicollar, the Seifert form on $H_{1}\left(
F\right) $ is defined to be the linking number of $x$ with $y^{+}$,
for any $x,y\in H_{1}\left( F\right) $, where $y^{+}$ denotes a
pushoff of $y$ in the positive direction of the bicollar of $F$.  A
Seifert matrix $V$ is the matrix representing the Seifert form with
respect to a choice of basis for $H_{1}\left( F\right) $.  For any
Seifert matrix $V$, recall that $V-tV^{T}$ presents the Alexander
module and $\det \left( V-V^{T}\right) =1\neq 0$. It follows that
the Alexander module is a torsion module.  That is, for any element
of the Alexander module, $x$, there is a non-zero element, $p\left(
t\right) $, of $ \mathbb{Z}\left[ t,t^{-1}\right] $ such that $x\ast
p\left( t\right) =0$.

Since the Alexander module is a torsion module, it is possible to
define a linking form on the Alexander module, known classically as
the Blanchfield linking form.  Let $x$ and $y$ be any two elements
of the Alexander module. Then as above, there is a non-zero element
of $\mathbb{Z}\left[ t,t^{-1} \right] $, $p\left( t\right) $, such
that $x\ast p\left( t\right) =0$.  Therefore there is a 2-chain,
$\alpha $, in the infinite cyclic cover of the exterior whose
boundary is $x\ast p\left( t\right) $.  We define the Blanchfield
linking form of $x$ and $y$ to be $\mathcal{B\ell }\left( x,y\right)
=\sum_{i=-\infty }^{\infty }\frac{1}{p\left( t^{-1}\right) } \lambda
\left( \alpha ,y\ast t^{i}\right) t^{-i}$ (mod $\mathbb{Z}\left[
t,t^{-1}\right] $), where $\lambda $ is the ordinary intersection
form. Notice that the Blanchfield linking form takes values in
$\mathbb{Q} \left( t\right) /\mathbb{Z}\left[ t,t^{-1}\right] $.

In order to motivate our main theorem, we recall some results about
the Blanchfield linking form. C. Kearton \cite{kearton} and H.F.
Trotter \cite{trotter} each proved the following theorem.

\begin{theorem}
If $V$ is a Seifert matrix for a knot $K$, then $\left( 1-t\right) \left[
V-tV^{T}\right] ^{-1}$ represents the Blanchfield linking form for $K$.
\end{theorem}

Recall that S-equivalence of matrices is the equivalence relation
generated by integral congruence and column enlargements.  Here $V$
is integrally congruent to $P^{T}VP$ where $P$ is an integral matrix
with $\det P=\pm 1$, and a column enlargement of $V$ is the
following.

\begin{equation*}
\left(
\begin{array}{ccc}
V & u^{T} & 0 \\
v & x & 1 \\
0 & 0 & 0
\end{array}
\right)
\end{equation*}
Here $x$ is an integer and $u$ and $v$ are column vectors.
Furthermore, two knots are S-equivalent if they have S-equivalent
Seifert matrices.

\begin{proposition}[\protect\cite{trotter}, p. 179; \protect\cite{kearton},
p. 142]
\label{Bl=Seq}Two knots have isomorphic Blanchfield linking forms if and
only if they are S-equivalent.
\end{proposition}

The question arises: Do there exist knots with isomorphic Alexander
modules, but non-isomorphic Blanchfield forms?  By Proposition
\ref{Bl=Seq}, it suffices to find examples of knots with isomorphic
Alexander modules that are not S-equivalent.  Furthermore, since the
ordinary signature of a knot is an S-equivalence invariant, we have
reduced the problem to finding two knots with isomorphic Alexander
modules, but with different signatures.

Given any knot, $K$, let $-K=r\overline{K}$ denote the reverse of
the mirror-image of $K$.  This is also the inverse of $K$ in the
knot concordance group.

\begin{proposition}
If $K$ is a knot such that the ordinary signature of $K$ is
non-zero, then $ K $ and $-K$ have isomorphic Alexander modules, but
non-isomorphic Blanchfield linking forms.
\end{proposition}

\begin{proof}
If $V$ is a Seifert matrix for $K$, then $-V^{T}$ is a Seifert
matrix for its mirror-image, $\overline{K}$, and $V^{T}$ is a
Seifert matrix for its reverse, $rK$.  Therefore the Seifert matrix
for $-K$ is $-V$.  Since $ V-tV^{T}$ and $-V+tV^{T}$ present
isomorphic modules, $K$ and $-K$ have isomorphic Alexander modules.
 However, if the signature of $K$ is non-zero, then the signature
of $-K$ is not equal to the signature of $K$.  Hence, $K$ and $-K$
are not S-equivalent, and therefore have non-isomorphic Blanchfield
linking forms.
\end{proof}

We note that the examples that Trotter provided in \cite{trotter} were found
using different methods than those presented here.

\section{Higher-Order Alexander Modules and Linking Forms}

Let us recall some of the definitions and results from \cite{nckt}
and \cite{cot1}.  Given a knot $K$, let $E\left( K\right) $ denote
the exterior of $ K $, $S^{3}\backslash \nbhd \left( K\right) $, and
let $G=\pi _{1}\left( E\left( K\right) \right) $.  Recall that the
derived series of a group $H$ is defined recursively by $H^{\left(
0\right) }=H$ and $H^{\left( n+1\right) }=\left[ H^{\left( n\right)
},H^{\left( n\right) }\right] $, for $ n\geq 1$.  We will use
$\Gamma _{n}$ to denote the quotient group $\frac{G}{ G^{\left(
n+1\right) }}$.  Then we have the coefficient system defined by the
homomorphism $G\rightarrow \Gamma _{n}$.

\begin{definition}
For $n\geq 0$, the $n$th higher-order Alexander module of a knot $K$
is
\begin{equation*}
\mathcal{A}_{n}\left( K\right) \equiv H_{1}\left( E\left( K\right)
;\mathbb{Z }\Gamma _{n}\right).
\end{equation*}
\end{definition}

Therefore the $n$th higher-order Alexander module is the first
(integral) homology group of the covering space of the knot exterior
corresponding to $ G^{\left( n+1\right) }$, considered as a right
$\mathbb{Z}\Gamma _{n}$-module.  This is the same as
$\frac{G^{\left( n+1\right) }}{G^{\left( n+2\right) }}$ as a right
$\mathbb{Z}\Gamma _{n}$-module.  (Notice that we are working with
the right $\mathbb{Z}\Gamma _{n}$-module structure on the chain
groups given by $\alpha \ast g=g\alpha g^{-1}$.)

As in the classical case, the higher-order Alexander modules of a
knot are torsion modules (\cite{nckt}, Prop. 3.10).  Therefore, it
is possible to define symmetric linking forms on the higher-order
Alexander modules.  Before giving the careful algebraic definition,
we describe the geometric idea of the linking forms.  Given any
element $x\in \mathcal{A}_{n}\left( K\right) $, there is some
$\gamma \in \mathbb{Z}\Gamma _{n}$ such that $ x\cdot \gamma =0$.
Therefore $x\cdot \gamma $ can be represented as the boundary of a
2-chain, $\alpha \in C_{2}\left( X;\mathbb{Z}\Gamma _{n}\right) $.
For any $y\in \mathcal{A}_{n}\left( K\right) $, define $
\mathcal{B\ell }_{n}\left( x,y\right) \equiv \overline{\gamma
}^{-1}\cdot \lambda _{n}\left( \alpha ,y\right) $, where $\lambda
_{n}$ denotes the equivariant intersection pairing on $E\left(
K\right) $ with coefficients in $\mathbb{Z}\Gamma _{n}$.  Here
$\overline{\gamma }$ is the image of $\gamma $ under the group ring
involution defined on the group ring $\mathbb{Z} \Gamma _{n}$ by
$\overline{\left( \sum n_{i}g_{i}\right) }=\sum n_{i}g_{i}^{-1}$
(see \cite[p.5]{passman}).  Since we are working with the right
module structure, the equivariant intersection pairing is defined as
$$ \lambda _{n}\left( \alpha ,y\right) =\sum\limits_{g\in \Gamma
_{n}}\lambda \left( \alpha ,\left( y\ast g\right) \right) \cdot
g^{-1},$$ where $\lambda $ denotes the ordinary intersection form.

In \cite[Prop. 3.2]{cot1}, it is shown that $\mathbb{Z}\Gamma _{n}$
is an Ore domain.  Therefore it is possible to define the right ring
of fractions of $\mathbb{Z}\Gamma _{n}$ (see \cite[Cor.
1.3.3]{cohn}), which we will denote by $\mathcal{K}_{n}$.  The short
exact sequence $0\rightarrow \mathbb{Z}\Gamma _{n}\rightarrow
\mathcal{K}_{n}\rightarrow \mathcal{K}_{n}/ \mathbb{Z}\Gamma
_{n}\rightarrow 0$ gives rise to the Bockstein sequence.
\begin{equation*}
H_{2}\left( E\left( K\right) ;\mathcal{K}_{n}\right) \rightarrow
H_{2}\left( E\left( K\right) ;\mathcal{K}_{n}/\mathbb{Z}\Gamma
_{n}\right) \overset{B}{ \rightarrow }H_{1}\left( E\left( K\right)
;\mathbb{Z}\Gamma _{n}\right) \rightarrow H_{1}\left( E\left(
K\right) ;\mathcal{K}_{n}\right)
\end{equation*}
Since the higher-order Alexander modules of a knot are torsion
modules, it follows that $H_{i}\left( E\left( K\right)
;\mathcal{K}_{n}\right) =0$, for $ i=1,2$ (\cite{nckt}, Cor. 3.12).
Therefore the Bockstein map, $ B:H_{2}\left( E\left( K\right)
;\mathcal{K}_{n}/\mathbb{Z}\Gamma _{n}\right) \rightarrow
H_{1}\left( E\left( K\right) ;\mathbb{Z}\Gamma _{n}\right) $, is an
isomorphism.

Let $\mathcal{A}_{n}\left( K\right) ^{\#}\,\equiv \overline{\Hom_{
\mathbb{Z}\Gamma _{n}}\left( \mathcal{A}_{n}\left( K\right)
,\mathcal{K}_{n}/ \mathbb{Z}\Gamma _{n}\right) }$, where given any
left $\mathcal{R}$-module $ \mathcal{M}$, $\overline{\mathcal{M}}$
represents the usual associated right $\mathcal{R}$-module resulting
from the involution of $\mathcal{R}$.  We now give the precise
definition of the symmetric linking forms defined on the
higher-order Alexander modules.

\begin{definition}
The $n$th higher-order linking form, $\mathcal{B\ell
}_{n}:\mathcal{A} _{n}\left( K\right) \rightarrow
\mathcal{A}_{n}\left( K\right) ^{\#}$, for a knot $K$, is the
composition of the following maps:
\begin{eqnarray*}
&&\mathcal{A}_{n}\left( K\right) \overset{B^{-1}}{\rightarrow
}H_{2}\left( E\left( K\right) ;\mathcal{K}_{n}/\mathbb{Z}\Gamma
_{n}\right) \overset{\pi } {\rightarrow }H_{2}\left( E\left(
K\right) ,\partial E\left( K\right) ;
\mathcal{K}_{n}/\mathbb{Z}\Gamma _{n}\right) \\
&&\hspace{0.75in}\overset{P.D.}{\rightarrow }\overline{H^{1}\left(
E\left( K\right) ;\mathcal{K}_{n}/\mathbb{Z}\Gamma _{n}\right)
}\overset{\kappa }{ \rightarrow }\mathcal{A}_{n}\left( K\right)
^{\#}
\end{eqnarray*}
where $P.D.$ is the Poincar\'{e} duality isomorphism and $\kappa $
is the Kronecker evaluation map. \ We will often denote $\left[
\mathcal{B\ell } _{n}\left( x\right) \right] \left( y\right) $ by
$\mathcal{B\ell }_{n}\left( x,y\right) $.
\end{definition}

The higher-order linking forms that we consider differ from those
defined by T. Cochran in \cite{nckt} and T. Cochran, K. Orr, and P.
Teichner in \cite{cot1} because we do not localize the coefficients
and because the higher-order linking forms that we consider are
canonically associated to $K$, unlike those that were the focus of
\cite{cot1}.  Furthermore, since $\mathbb{ Z}\Gamma _{n}$ is not a
PID, our linking forms may be singular.

\section{Genetic Infection}

In order to construct the desired examples, we use a satellite
technique, that was called \textit{genetic infection} in
\cite{nckt}.  Let $K$ and $J$ be fixed knots, and let $\eta $ be an
embedded oriented circle in $ S^{3}\backslash K$ which is itself
unknotted in $S^{3}$.  Since $\eta $ is unknotted in $S^{3}$, it
bounds a disc, $D$, in $S^{3}$, which we can choose to intersect $K$
transversely. We construct a new knot by tying the strands of $K$
that pierce $D$ into the knot $J$.  That is, we replace the strands
of $K$ that intersect a small neighborhood of $D$ with untwisted
parallels of a knotted arc with oriented knot type $J$.  We call the
resulting knot \textit{the result of infecting }$K$ \textit{by }$J$
\textit{along }$\eta $, denoted by $K\left( \eta ,J\right) $.
Alternatively, we can view this construction from a surgery point of
view.  Beginning with the exterior of $ K$, $E\left( K\right) $,
delete the interior of a tubular neighborhood of $ \eta $, and
replace it with the exterior of $J$, $E\left( J\right) $,
identifying the meridian of $J$ with the longitude of $\eta $, and
the longitude of $J$ with the inverse of the meridian of $\eta $.
The result is the exterior of $K\left( \eta ,J\right) $.  This
surgery description is better suited for our purposes.

Since there is a degree one map (rel boundary) $E\left( J\right) \rightarrow
E\left( \text{unknot}\right) $, there is a degree one map $f:E\left( K\left(
\eta ,J\right) \right) \rightarrow E\left( K\right) $, which is the identity
outside of $E\left( J\right) $.

\begin{proposition}[\protect\cite{nckt}, Thm. 8.1]
\label{nckt8.1}If $\eta \in \pi _{1}\left( E\left( K\right) \right)
^{\left( n\right) }$, then the map $f$ induces an isomorphism:
\begin{equation*}
f_{\ast }:\frac{\pi _{1}\left( E\left( K\left( \eta ,J\right)
\right) \right) }{\pi _{1}\left( E\left( K\left( \eta ,J\right)
\right) \right) ^{\left( n+1\right) }}\rightarrow \frac{\pi
_{1}\left( E\left( K\right) \right) }{\pi _{1}\left( E\left(
K\right) \right) ^{\left( n+1\right) }}.
\end{equation*}
\end{proposition}

Therefore we will use $\Gamma _{n}$ to denote both groups.  The
following composition of maps defines coefficient systems on
$E\left( J\right) $, $ E\left( K\left( \eta ,J\right) \right) $, and
$E\left( K\right) $.
\begin{equation*}
\pi _{1}\left( E\left( J\right) \right) \overset{i_{\ast
}}{\rightarrow }\pi _{1}\left( E\left( K\left( \eta ,J\right)
\right) \right) \overset{f_{\ast }} {\rightarrow }\pi _{1}\left(
E\left( K\right) \right) \overset{\phi }{ \twoheadrightarrow }\Gamma
_{n}
\end{equation*}

The following results demonstrate the relationship between genetic
infection and the higher-order Alexander modules.  We include the
proofs from \cite{nckt} because some of our results will be proved
using similar techniques.

\begin{corollary}
\label{nckt8.1b}If $\eta \in \pi _{1}\left( E\left( K\right) \right)
^{\left( n\right) }$, then $f:E\left( K\left( \eta ,J\right) \right)
\rightarrow E\left( K\right) $ induces isomorphisms between the
$i$-th order Alexander modules of $K\left( \eta ,J\right) $ and $K$,
for $0\leq i\leq n-1$.
\end{corollary}

\begin{proof}
Let $G=\pi _{1}\left( E\left( K\right) \right) $ and
$\widehat{G}=\pi _{1}\left( E\left( K\left( \eta ,J\right) \right)
\right) $. \ We have the following commutative diagram with exact
rows.
\begin{equation*}
\begin{diagram} 1 & \rTo & \frac{\widehat{G}^{\left( i+1\right)
}}{\widehat{G}^{\left(i+2\right) }} & \rTo &
\frac{\widehat{G}}{\widehat{G}^{\left(i+2\right) }} & \rTo &
\frac{\widehat{G}}{\widehat{G}^{\left(i+1\right) }} & \rTo & 1 \\ & &
\dTo>{f_{\ast}} & & \dTo>{f_{\ast}}<{\cong} & & \dTo>{f_{\ast}}<{\cong} & &
\\ 1 & \rTo & \frac{G^{\left( i+1\right) }}{G^{\left(i+2\right) }} & \rTo &
\frac{G}{G^{\left(i+2\right) }} & \rTo & \frac{G}{G^{\left(i+1\right) }} &
\rTo & 1 \\ \end{diagram}
\end{equation*}
For $0\leq i\leq n-1$, the middle and right vertical maps are
isomorphisms by Proposition \ref{nckt8.1}. \ Therefore $f_{\ast
}:\frac{\widehat{G} ^{\left( i+1\right) }}{\widehat{G}^{\left(
i+2\right) }}\rightarrow \frac{ G^{\left( i+1\right) }}{G^{\left(
i+2\right) }}$ is an isomorphism.  That is, $\mathcal{A}_{i}\left(
K\left( \eta ,J\right) \right) \cong \mathcal{A} _{i}\left( K\right)
$.
\end{proof}

In the proof of the next theorem we will require following lemma,
which we state without proof.

\begin{lemma}[\protect\cite{nckt}, Lemma 8.3]
\label{nckt8.3}If $\eta \in \pi _{1}\left( E\left( K\right) \right)
^{\left( n\right) }$ and $\eta \notin \pi _{1}\left( E\left(
K\right) \right) ^{\left( n+1\right) }$, then the inclusion
$\partial E\left( J\right) \rightarrow E\left( J\right) $ induces an
isomorphism on $H_{0}\left( -; \mathbb{Z}\Gamma _{n}\right) $ and
the trivial map\ on $H_{1}\left( -; \mathbb{Z}\Gamma _{n}\right) $.
\end{lemma}

\begin{theorem}[\protect\cite{nckt}, Thm. 8.2]
\label{nckt8.2}If $\eta \in \pi _{1}\left( E\left( K\right) \right)
^{\left( n\right) }$ and $\eta \notin \pi _{1}\left( E\left(
K\right) \right) ^{\left( n+1\right) }$, then
\begin{equation*}
H_{1}\left( E\left( K\left( \eta ,J\right) \right) ;\mathbb{Z}\Gamma
_{n}\right) \cong \mathcal{A}_{n}\left( K\right) \oplus H_{1}\left( E\left(
J\right) ;\mathbb{Z}\Gamma _{n}\right) \text{.}
\end{equation*}
\end{theorem}

\begin{proof}

Let $E\left( \eta \right) $ denote the result of deleting the
interior of a tubular neighborhood of $\eta $ from the exterior of
$K$. Using the surgery description of genetic infection, we have
$E\left( K\left( \eta ,J\right) \right) \cong E\left( J\right) \cup
_{\partial E\left( J\right) }E\left( \eta \right) $.  Since
infecting with the unknot leaves the knot unchanged, we can view the
exterior of $K$ as the union of the exterior of the unknot, $ U$,
and $E\left( \eta \right) $.  That is, $E\left( K\left( \eta
,J\right) \right) \cong E\left( U\right) \cup _{\partial E\left(
J\right) }E\left( \eta \right) $.  Of course, $E\left( U\right) $ is
just a solid torus.

We consider the Mayer-Vietoris sequence with $\mathbb{Z}\Gamma _{n}$
coefficients for $E\left( K\left( \eta ,J\right) \right) \cong
E\left( J\right) \cup _{\partial E\left( J\right) }E\left( \eta
\right) $.
\begin{equation*}
\begin{diagram} H_{1}\left( E\left( K\left(\eta,J\right)\right) \right) &
\rTo^{\partial _{\ast }} & H_{0}\left( \partial E\left( J\right) \right) &
\rTo^{\left( \psi _{1},\psi _{2}\right) } & H_{0}\left( E\left( J\right)
\right) \oplus H_{0}\left( E\left( \eta \right) \right) \\ \end{diagram}
\end{equation*}
By Lemma \ref{nckt8.3}, $\psi _{1}:H_{0}\left( \partial E\left(
J\right) \right) \rightarrow H_{0}\left( E\left( J\right) \right) $
is an isomorphism.  Therefore $\partial _{\ast }:H_{1}\left( E\left(
K\left( \eta ,J\right) \right) \right) \rightarrow H_{0}\left(
\partial E\left( J\right) \right) $ is the trivial map.  Similarly,
since $E\left( K\right) \cong E\left( U\right) \cup _{\partial
E\left( J\right) }E\left( \eta \right) $, $
\partial _{\ast }:H_{1}\left( E\left( K\right) \right) \rightarrow
H_{0}\left( \partial E\left( J\right) \right) $ is the trivial map.
Thus we have the following diagram.
\begin{equation*}
\begin{diagram} H_{1}\left( \partial E\left( J\right) \right) & \rTo^{\left(
\psi _{1},\psi _{2}\right) } & H_{1}\left( E\left( J\right) \right) \oplus
H_{1}\left( E\left( \eta \right) \right) & \rTo^{i_{\ast }+j_{\ast }} &
H_{1}\left( E\left(K\left(\eta,J\right)\right) \right) & \rTo^{\partial
_{\ast }} & 0 \\ \dTo>{f_{\ast}} & & \dTo>{f_{\ast}} & & \dTo>{f_{\ast}} & &
\\ H_{1}\left( \partial E\left( J\right) \right) & \rTo^{\left( \psi
_{1}^{\prime },\psi _{2}\right) } & H_{1}\left( E\left( U\right) \right)
\oplus H_{1}\left( E\left( \eta \right) \right) & \rTo^{i_{\ast
}^{\prime}+j_{\ast }} & H_{1}\left( E\left( K\right) \right) &
\rTo^{\partial _{\ast }} & 0 \\ \end{diagram}
\end{equation*}

Lemma \ref{nckt8.3} states that $\psi _{1}:H_{1}\left( \partial
E\left( J\right) \right) \rightarrow H_{1}\left( E\left( J\right)
\right) $ is the trivial map.  Therefore $\im\left( \psi _{1},\psi
_{2}\right) =0\oplus \psi _{2}\left( H_{1}\left( \partial E\left(
J\right) \right) \right) $.  Hence $H_{1}\left( E\left( K\left( \eta
,J\right) \right) \right) \cong H_{1}\left( E\left( J\right) \right)
\oplus \frac{H_{1}\left( E\left( \eta \right) \right) }{\psi
_{2}\left( H_{1}\left( \partial E\left( J\right) \right) \right) }$.
 Similarly, $H_{1}\left( E\left( K\right) \right) \cong H_{1}\left(
E\left( U\right) \right) \oplus \frac{H_{1}\left( E\left( \eta
\right) \right) }{\psi _{2}\left( H_{1}\left( \partial E\left(
J\right) \right) \right) }$.  Notice that since $f$ is the identity
on $\partial E\left( J\right) $ and $E\left( \eta \right) $,
$f_{\ast }\circ \psi _{2}=\psi _{2}\circ f_{\ast }$.  However, since
$\eta \notin \pi _{1}\left( E\left( K\right) \right) ^{\left(
n+1\right) }$, $\mu _{J}$, which is the generator of $\pi _{1}\left(
E\left( U\right) ;\mathbb{Z} \right) $, gets unwound in the $\Gamma
_{n}$-cover. Hence $H_{1}\left( E\left( U\right) \right) =0$.
Therefore $H_{1}\left( E\left( K\left( \eta ,J\right) \right)
;\mathbb{Z}\Gamma _{n}\right) \cong H_{1}\left( E\left( J\right)
;\mathbb{Z}\Gamma _{n}\right) \oplus H_{1}\left( E\left( K\right) ;
\mathbb{Z}\Gamma _{n}\right) .$
\end{proof}

\begin{corollary}
\label{SplitSES}If $\eta \in \pi _{1}\left( E\left( K\right) \right)
^{\left( n\right) }$ and $\eta \notin \pi _{1}\left( E\left( K\right)
\right) ^{\left( n+1\right) }$, then
\begin{equation*}
0\rightarrow H_{1}\left( E\left( J\right) ;\mathbb{Z}\Gamma
_{n}\right) \overset{i_{\ast }}{\rightarrow }H_{1}\left( E\left(
K\left( \eta ,J\right) \right) ;\mathbb{Z}\Gamma _{n}\right)
\overset{f_{\ast }}{\rightarrow } H_{1}\left( E\left( K\right)
;\mathbb{Z}\Gamma _{n}\right) \rightarrow 0
\end{equation*}
is a split short exact sequence.
\end{corollary}

Since $\pi _{1}\left( E\left( J\right) \right) $ is normally
generated by the meridian of $J$, it follows that if $\eta \in \pi
_{1}\left( E\left( K\right) \right) ^{\left( n\right) }$, then the
image of $\pi _{1}\left( E\left( J\right) \right) $ is contained in
$\pi _{1}\left( E\left( K\right) \right) ^{\left( n\right) }$.
Therefore, $\pi _{1}\left( E\left( J\right) \right) ^{\prime }$ is
in $\pi _{1}\left( E\left( K\right) \right) ^{\left( n+1\right) }$,
and thus in the kernel of the composition: $$\pi _{1}\left( E\left(
J\right) \right) \overset{i_{\ast }}{\rightarrow }\pi _{1}\left(
E\left( K\left( \eta ,J\right) \right) \right) \overset{f_{\ast }}{
\rightarrow }\pi _{1}\left( E\left( K\right) \right) \overset{\phi
}{ \twoheadrightarrow }\Gamma _{n}.$$  Hence, we have a ring
homomorphism: $$ \mathbb{Z}\left[ \frac{\pi _{1}\left( E\left(
J\right) \right) }{\pi _{1}\left( E\left( J\right) \right) ^{\prime
}}\right] \cong \mathbb{Z}\left[ t,t^{-1}\right] \overset{\psi
}{\rightarrow }\mathbb{Z}\Gamma _{n}.$$  If $ \eta \notin \pi
_{1}\left( E\left( K\right) \right) ^{n+1}$, this is a monomorphism.

\begin{corollary}
\label{nckt8.2TimVersion}If $\eta \in \pi _{1}\left( E\left(
K\right) \right) ^{\left( n\right) }$, then $\mathcal{A}_{i}\left(
K\left( \eta ,J\right) \right) \cong \mathcal{A}_{i}\left( K\right)
$
\begin{equation*}
\mathcal{A}_{n}\left( K\left( \eta ,J\right) \right) \cong
\mathcal{A} _{n}\left( K\right) \oplus \left( \mathcal{A}_{0}\left(
J\right) \otimes _{ \mathbb{Z}\left[ t,t^{-1}\right]
}\mathbb{Z}\Gamma _{n}\right)
\end{equation*}
where $\mathbb{Z}\Gamma _{n}$ is a left $\mathbb{Z}\left[
t,t^{-1}\right]$-module by the homomorphism sending $t$ to $\phi
\left( \eta \right) $.
\end{corollary}

\begin{proof}
If $\eta \in \pi _{1}\left( E\left( K\right) \right) ^{\left(
n+1\right) }$, then $\mathcal{A}_{n}\left( K\left( \eta ,J\right)
\right) \cong \mathcal{A} _{n}\left( K\right) $ by Corollary
\ref{nckt8.1b}.  Also in this case, $ \phi \left( \eta \right) =1$.
 Therefore, $\mathcal{A}_{0}\left( J\right)
\otimes _{\mathbb{Z}\left[ t,t^{-1}\right] }\mathbb{Z}\Gamma
_{n}=0$, since $ \mathcal{A}_{0}\left( J\right) $ is presented by
$V-tV^{t}$ and $\det \left( V-V^{t}\right) =1$, where $V$ is a
Seifert matrix for $J$.

If $\eta \notin \pi _{1}\left( E\left( K\right) \right) ^{\left(
n+1\right) } $, then $H_{1}\left( E\left( K\left( \eta ,J\right)
\right) ;\mathbb{Z} \Gamma _{n}\right) \cong \mathcal{A}_{n}\left(
K\right) \oplus H_{1}\left( E\left( J\right) ;\mathbb{Z}\Gamma
_{n}\right) $ by Theorem \ref{nckt8.2}.
 Furthermore, if $\widetilde{E\left( J\right) }$ is the universal cover of $
E\left( J\right) $,
\begin{eqnarray*}
H_{1}\left( E\left( J\right) ;\mathbb{Z}\Gamma _{n}\right)
&=&H_{1}\left( C_{\ast }\left( \widetilde{E\left( J\right) }\right)
\otimes _{\mathbb{Z}\pi
_{1}\left( E\left( J\right) \right) }\mathbb{Z}\Gamma _{n}\right) \\
&\cong& H_{1}\left( C_{\ast }\left( \widetilde{E\left( J\right)
}\right) \otimes _{\mathbb{Z}\pi _{1}\left( E\left( J\right) \right)
}\mathbb{Z}\left[ t,t^{-1}\right] \otimes _{\mathbb{Z}\left[
t,t^{-1}\right] }\mathbb{Z}\Gamma _{n}\right) \\
&&\hspace{2in}\text{ \cite[p. 109]{hilton-stammbach}}
\end{eqnarray*}
Since $\eta \notin \pi _{1}\left( E\left( K\right) \right) ^{n+1}$,
it follows that $\mathbb{Z}\left[ t,t^{-1}\right] \overset{\psi
}{\rightarrow } \mathbb{Z}\Gamma _{n}$ is a monomorphism.  It
follows from \cite[Lemma 1.3]{passman} that $\mathbb{Z}\Gamma _{n}$
is a free, and therefore flat $ \mathbb{Z}\left[ t,t^{-1}\right]
$-module.  Hence,
\begin{eqnarray*}
H_{1}\left( E\left( J\right) ;\mathbb{Z}\Gamma _{n}\right) &\cong
&H_{1}\left( C_{\ast }\left( \widetilde{E\left( J\right) }\right)
\otimes _{ \mathbb{Z}\pi _{1}\left( E\left( J\right) \right)
}\mathbb{Z}\left[ t,t^{-1} \right] \right) \otimes
_{\mathbb{Z}\left[ t,t^{-1}\right] }\mathbb{Z}\Gamma
_{n} \\
&\cong &H_{1}\left( E\left( J\right) ;\mathbb{Z}\left[ t,t^{-1}\right]
\right) \otimes _{\mathbb{Z}\left[ t,t^{-1}\right] }\mathbb{Z}\Gamma _{n} \\
&\cong &\mathcal{A}_{0}\left( J\right) \otimes _{\mathbb{Z}\left[
t,t^{-1} \right] }\mathbb{Z}\Gamma _{n}
\end{eqnarray*}
\end{proof}

\section{The Effect of Genetic Infection on the Higher-Order Linking Forms
\label{SectionMarker}}

The idea behind the construction of our examples is to infect the
same knot $ K$ along the same element $\eta \in \pi _{1}\left(
E\left( K\right) \right) ^{\left( n\right) }$, $\eta \notin \pi
_{1}\left( E\left( K\right) \right) ^{\left( n+1\right) }$ by two
different knots $J_{1}$, $J_{2}$ that have isomorphic classical
Alexander modules.  Corollary \ref{nckt8.2TimVersion} implies that
the results of these infections will have isomorphic $i$-th order
Alexander modules for $0\leq i\leq n$.  We need to choose the knots
so that the higher-order linking forms of the results of the
infections are not isomorphic.  In this section, we will determine
the effect of genetic infection on the higher-order linking forms in
order to determine what the desired conditions are on the infecting
knots.

\begin{theorem}
\label{lowerlinkingforms}If $\eta \in \pi _{1}\left( E\left( K\right)
\right) ^{\left( n\right) }$, then $f:E\left( K\left( \eta ,J\right) \right)
\rightarrow E\left( K\right) $ induces isomorphisms between the $i$-th order
linking forms of $K\left( \eta ,J\right) $ and $K$, for $0\leq i\leq n-1$.
\end{theorem}

\begin{proof}
We have the following diagram.
\begin{equation*}
\begin{diagram}[small] H_{1}\left( E\left( K\left( \eta ,J\right) \right)
;\mathbb{Z}\Gamma _{i}\right) & \rTo_{f_{\ast}}^{\cong} & H_{1}\left(
E\left( K\right) ;\mathbb{Z}\Gamma _{i}\right) \\ \dTo>{B^{-1}} & &
\dTo>{B^{-1}} \\ H_{2}\left( E\left( K\left( \eta ,J\right) \right)
;\mathcal{K}_{i}/\mathbb{Z}\Gamma _{i}\right) & \rTo_{f_{\ast}} &
H_{2}\left( E\left( K\right) ;\mathcal{K}_{i}/\mathbb{Z}\Gamma _{i}\right)
\\ \dTo>{\pi} & & \dTo>{\pi} \\ H_{2}\left( E\left( K\left( \eta ,J\right)
\right) ,\partial E\left( K\left( \eta ,J\right) \right)
;\mathcal{K}_{i}/\mathbb{Z}\Gamma _{i}\right) & \rTo_{f_{\ast}} &
H_{2}\left( E\left( K\right) ,\partial E\left( K\right)
;\mathcal{K}_{i}/\mathbb{Z}\Gamma _{i}\right) \\ \dTo>{P.D.} & & \dTo>{P.D.}
\\ \overline{H^{1}\left( E\left( K\left( \eta ,J\right) \right)
;\mathcal{K}_{i}/\mathbb{Z}\Gamma _{i}\right) } & \lTo_{f^{\ast}} &
\overline{H^{1}\left( E\left( K\right) ;\mathcal{K}_{i}/\mathbb{Z}\Gamma
_{i}\right)} \\ \dTo>{\kappa} & & \dTo>{\kappa} \\ H_{1}\left( E\left(
K\left( \eta ,J\right) \right) ;\mathbb{Z}\Gamma _{i}\right) ^{\#} &
\lTo_{f^{\ast}}^{\cong} & H_{1}\left( E\left( K\right) ;\mathbb{Z}\Gamma
_{i}\right) ^{\#} \\ \end{diagram}
\end{equation*}
By the naturality of the Bockstein isomorphism, $\pi $, the Poincar\'{e}
duality isomorphism, and the Kronecker map, $\mathcal{B\ell }_{i}\left(
K\left( \eta ,J\right) \right) =f^{\ast }\circ \mathcal{B\ell }_{i}\left(
K\right) \circ f_{\ast }$. \ Since, by Corollary \ref{nckt8.1b}, $f$ induces
isomorphisms between the $i$-th order Alexander modules of $K\left( \eta
,J\right) $ and $K$, for $0\leq i<n$, it follows that $f$ induces
isomorphisms between $\mathcal{B\ell }_{i}\left( K\left( \eta ,J\right)
\right) $ and $\mathcal{B\ell }_{i}\left( K\right) $.
\end{proof}

From now on, we will regard $n$ as fixed and restrict our attention
to the case where $\eta \in \pi _{1}\left( E\left( K\right) \right)
^{\left( n\right) }$ and $\eta \notin \pi _{1}\left( E\left(
K\right) \right) ^{\left( n+1\right) }$.  As a result, we will
suppress $n$ from our notation for the higher-order linking forms.
Therefore, let $\mathcal{ B\ell }_{K\left( \eta ,J\right)
}:H_{1}\left( E\left( K\left( \eta ,J\right) \right)
;\mathbb{Z}\Gamma _{n}\right) \rightarrow H_{1}\left( E\left(
K\left( \eta ,J\right) \right) ;\mathbb{Z}\Gamma _{n}\right) ^{\#}$
and $ \mathcal{B\ell }_{K}:H_{1}\left( E\left( K\right)
;\mathbb{Z}\Gamma _{n}\right) \rightarrow H_{1}\left( E\left(
K\right) ;\mathbb{Z}\Gamma _{n}\right) ^{\#}$ denote the
higher-order linking forms for $K\left( \eta ,J\right) $ and $K$,
respectively.

Recall that the following composition of maps defines a coefficient
system on $E\left( J\right) $.
\begin{equation*}
\pi _{1}\left( E\left( J\right) \right) \overset{i_{\ast
}}{\rightarrow }\pi _{1}\left( E\left( K\left( \eta ,J\right)
\right) \right) \overset{f_{\ast }} {\rightarrow }\pi _{1}\left(
E\left( K\right) \right) \overset{\phi }{ \twoheadrightarrow }\Gamma
_{n}
\end{equation*}
This coefficient system is non-trivial if $\eta \notin \pi
_{1}\left( E\left( K\right) \right) ^{\left( n+1\right) }$.

\begin{proposition}
\label{ZGnTorsion}If $\eta \in \pi _{1}\left( E\left( K\right)
\right) ^{\left( n\right) }$ and $\eta \notin \pi _{1}\left( E\left(
K\right) \right) ^{\left( n+1\right) }$, then $H_{1}\left( E\left(
J\right) ;\mathbb{Z }\Gamma _{n}\right) $ is a $\mathbb{Z}\Gamma
_{n}$-torsion module.
\end{proposition}

\begin{proof}
Since $\eta \in \pi _{1}\left( E\left( K\right) \right) ^{\left(
n\right) }$ and $\eta \notin \pi _{1}\left( E\left( K\right) \right)
^{\left( n+1\right) }$ we have a ring monomorphism $\mathbb{Z}\left[
t,t^{-1}\right] \rightarrow \mathbb{Z}\Gamma _{n}$.  Therefore we
have $\mathbb{Z}\left[ t,t^{-1}\right] $-module homomorphisms
$\mathbb{Q}\left( t\right) \hookrightarrow \mathcal{K}_{n}$.
Furthermore, from the proof of Corollary \ref{nckt8.2TimVersion}, $
H_{1}\left( E\left( J\right) ;\mathbb{Z}\Gamma _{n}\right) \cong
\mathcal{A} _{0}\left( J\right) \otimes _{\mathbb{Z}\left[
t,t^{-1}\right] }\mathbb{Z} \Gamma _{n}$ where $\mathbb{Z}\Gamma
_{n}$ is a left $\mathbb{Z}\left[ t,t^{-1}\right] $-module by the
homomorphism sending $t$ to $\phi \left( \eta \right) $.  Hence, we
have the following.
\begin{eqnarray*}
\mathcal{A}_{0}\left( J\right) \otimes _{\mathbb{Z}\left[
t,t^{-1}\right] } \mathbb{Z}\Gamma _{n}\otimes _{\mathbb{Z}\Gamma
_{n}}\mathcal{K}_{n} &\cong & \mathcal{A}_{0}\left( J\right) \otimes
_{\mathbb{Z}\left[ t,t^{-1}\right] }
\mathcal{K}_{n} \\
&\cong &\mathcal{A}_{0}\left( J\right) \otimes _{\mathbb{Z}\left[
t,t^{-1} \right] }\mathbb{Q}\left( t\right) \otimes
_{\mathbb{Q}\left( t\right) }
\mathcal{K}_{n} \\
&\cong &0\text{ \ since }\mathcal{A}_{0}\left( J\right) \text{ is a
}\mathbb{ Z}\left[ t,t^{-1}\right] \text{-torsion module}
\end{eqnarray*}
\end{proof}

Again, the short exact sequence $0\rightarrow \mathbb{Z}\Gamma
_{n}\rightarrow \mathcal{K}_{n}\rightarrow
\mathcal{K}_{n}/\mathbb{Z}\Gamma _{n}\rightarrow 0$ gives rise to a
Bockstein sequence.
\begin{equation*}
H_{2}\left( E\left( J\right) ;\mathcal{K}_{n}\right) \rightarrow
H_{2}\left( E\left( J\right) ;\mathcal{K}_{n}/\mathbb{Z}\Gamma
_{n}\right) \overset{B}{ \rightarrow }H_{1}\left( E\left( J\right)
;\mathbb{Z}\Gamma _{n}\right) \rightarrow H_{1}\left( E\left(
J\right) ;\mathcal{K}_{n}\right)
\end{equation*}

\begin{corollary}
If $\eta \in \pi _{1}\left( E\left( K\right) \right) ^{\left(
n\right) }$ and $\eta \notin \pi _{1}\left( E\left( K\right) \right)
^{\left( n+1\right) }$, then for $i=1,2$, $$H_{i}\left( E\left(
J\right) ;\mathcal{K}_{n}\right) =0.$$
\end{corollary}

\begin{proof}
Since $\mathcal{K}_{n}$ is a flat $\mathbb{Z}\Gamma _{n}$-module
\cite[Prop. II.3.5]{stenstrom}, $H_{1}\left( E\left( K\right)
;\mathcal{K}_{n}\right) \cong \linebreak H_{1}\left( E\left(
K\right) ;\mathbb{Z}\Gamma _{n}\right) \otimes _{\mathbb{Z}\Gamma
_{n}}\mathcal{K}_{n}=0$ since $\mathcal{A} _{n}\left( K\right) $ is
a torsion module.  Similarly, it follows from Prop. 3.7 of
\cite{nckt} that $H_{0}\left( \partial E\left( K\right) ;
\mathcal{K}_{n}\right) =0$.  Therefore by the long exact sequence of
a pair, $H_{1}\left( E\left( K\right) ,\partial E\left( K\right)
;\mathcal{K} _{n}\right) =0$.  By Poincar\'{e} duality and the
Universal Coefficient Theorem for modules over the (noncommutative)
principal ideal domain $ \mathcal{K}_{n}$ \cite[pp.
44,102]{daviskirk}, we have: $$H_{2}\left( E\left( K\right)
;\mathcal{K}_{n}\right) \cong \Hom_{\mathcal{K}_{n}}\left(
H_{1}\left( E\left( K\right) ,\partial E\left( K\right) ;\mathcal{K}
_{n}\right) ,\mathcal{K}_{n}\right) =0.$$\end{proof}

\begin{corollary}
If $\eta \in \pi _{1}\left( E\left( K\right) \right) ^{\left(
n\right) }$ and $\eta \notin \pi _{1}\left( E\left( K\right) \right)
^{\left( n+1\right) }$, then the Bockstein map, $B:H_{2}\left(
E\left( J\right) ;\mathcal{K} _{n}/\mathbb{Z}\Gamma _{n}\right)
\rightarrow H_{1}\left( E\left( J\right) ; \mathbb{Z}\Gamma
_{n}\right) $, is an isomorphism.
\end{corollary}

\begin{definition}
If $\eta \in \pi _{1}\left( E\left( K\right) \right) ^{\left(
n\right) }$ and $\eta \notin \pi _{1}\left( E\left( K\right) \right)
^{\left( n+1\right) }$, define $\mathcal{B\ell }_{K\left( \eta
,J\right) }^{\otimes }:H_{1}\left( E\left( J\right)
;\mathbb{Z}\Gamma _{n}\right) \rightarrow H_{1}\left( E\left(
J\right) ;\mathbb{Z}\Gamma _{n}\right) ^{\#}$, to be the composition
of the following maps:
\begin{eqnarray*}
&&H_{1}\left( E\left( J\right) ;\mathbb{Z}\Gamma _{n}\right)
\overset{B^{-1}} {\rightarrow }H_{2}\left( E\left( J\right)
;\mathcal{K}_{n}/\mathbb{Z}\Gamma _{n}\right) \overset{\pi
}{\rightarrow }H_{2}\left( E\left( J\right)
,\partial E\left( J\right) ;\mathcal{K}_{n}/\mathbb{Z}\Gamma _{n}\right) \\
&&\hspace{0.5in}\overset{P.D.}{\rightarrow }\overline{H^{1}\left(
E\left( J\right) ;\mathcal{K}_{n}/\mathbb{Z}\Gamma _{n}\right)
}\overset{\kappa }{ \rightarrow }H_{1}\left( E\left( J\right)
;\mathbb{Z}\Gamma _{n}\right) ^{\#}
\end{eqnarray*}
where $P.D.$ is the Poincar\'{e} duality isomorphism and $\kappa $
is the Kronecker evaluation map. We remark that the coefficient
system that we are using is defined using $K(\eta,J)$, and therefore
$\mathcal{B\ell }_{K\left( \eta ,J\right) }^{\otimes }$ does indeed
depend on $K$ and $\eta$, as well as $J$.
\end{definition}

Let $g$ be a splitting for the exact sequence in Corollary
\ref{SplitSES}. That is, $f_{\ast }\circ g=id$.

\begin{theorem}
\label{DirectSum}If $\eta \in \pi _{1}\left( E\left( K\right)
\right) ^{\left( n\right) }$ and $\eta \notin \pi _{1}\left( E\left(
K\right) \right) ^{\left( n+1\right) }$, then $\mathcal{B\ell
}_{K\left( \eta ,J\right) }\cong \mathcal{B\ell }_{K\left( \eta
,J\right) }^{\otimes }$ $\oplus \mathcal{B\ell }_{K}$.  That is,
\begin{equation*}
\mathcal{B\ell }_{K\left( \eta ,J\right) }^{\otimes }\left(
x_{1},y_{1}\right) +\mathcal{B\ell }_{K}\left( x_{2},y_{2}\right)
=\mathcal{ B\ell }_{K\left( \eta ,J\right) }\left( i_{\ast }\left(
x_{1}\right) +g\left( x_{2}\right) ,i_{\ast }\left( y_{1}\right)
+g\left( y_{2}\right) \right) \text{.}
\end{equation*}
\end{theorem}

\begin{proof}
We have the following diagram.
\begin{equation*}
\begin{diagram} H_{1}\left( E\left( J\right) ;\mathbb{Z}\Gamma _{n}\right) &
\rTo^{i_{\ast }} & H_{1}\left( E\left( K\left( \eta ,J\right) \right)
;\mathbb{Z}\Gamma _{n}\right) & \lTo^g & H_{1}\left( E\left( K\right)
;\mathbb{Z}\Gamma _{n}\right) \\ \dTo>{\mathcal{B\ell }_{K\left( \eta
,J\right) }^{\otimes }} & & \dTo>{\mathcal{B\ell }_{K\left( \eta ,J\right)
}} & & \dTo>{\mathcal{B\ell }_{K}} \\ H_{1}\left( E\left( J\right)
;\mathbb{Z}\Gamma _{n}\right)^{\#} & \lTo^{i^{\ast }} & H_{1}\left( E\left(
K\left( \eta ,J\right) \right) ;\mathbb{Z}\Gamma _{n}\right)^{\#} &
\rTo^{g^{\#}} & H_{1}\left( E\left( K\right) ;\mathbb{Z}\Gamma
_{n}\right)^{\#} \\ \end{diagram}
\end{equation*}
where $g^{\#}$ is the dual of $g$.  Notice that since $f_{\ast
}\circ g=id$, it follows that $g^{\#}\circ f^{\ast }=id$. \ The
isomorphism in the theorem will be given by $i_{\ast }\oplus
g:H_{1}\left( E\left( K\right) ;\mathbb{Z}\Gamma _{n}\right) \oplus
H_{1}\left( E\left( J\right) ;\mathbb{Z} \Gamma _{n}\right)
\rightarrow H_{1}\left( E\left( K\left( \eta ,J\right) \right)
;\mathbb{Z}\Gamma _{n}\right) $. \ Hence the theorem will follow
from the following four claims.

\begin{enumerate}
\item $g^{\#}\circ \mathcal{B\ell }_{K\left( \eta ,J\right) }\circ g=
\mathcal{B\ell }_{K}$ \newline which establishes $\mathcal{B\ell
}_{K\left( \eta ,J\right) }\left( g\left( x_{2}\right) ,g\left(
y_{2}\right) \right) = \mathcal{B\ell }_{K}\left( x_{2},y_{2}\right)
$

\item $i^{\ast }\circ \mathcal{B\ell }_{K\left( \eta ,J\right) }\circ
i_{\ast }=\mathcal{B\ell }_{K\left( \eta ,J\right) }^{\otimes }$
\newline which establishes $\mathcal{B\ell }_{K\left( \eta ,J\right)
}\left( i_{\ast }\left( x_{1}\right) ,i_{\ast }\left( y_{1}\right)
\right) =\mathcal{B\ell } _{K\left( \eta ,J\right) }^{\otimes
}\left( x_{1},y_{1}\right) $

\item $g^{\#}\circ \mathcal{B\ell }_{K\left( \eta ,J\right) }\circ i_{\ast
}=0$ which establishes $\mathcal{B\ell }_{K\left( \eta ,J\right) }\left(
i_{\ast }\left( x_{1}\right) ,g\left( y_{2}\right) \right) =0$

\item $i^{\ast }\circ \mathcal{B\ell }_{K\left( \eta ,J\right) }\circ g=0$
which establishes $\mathcal{B\ell }_{K\left( \eta ,J\right) }\left( g\left(
x_{2}\right) ,i_{\ast }\left( y_{1}\right) \right) =0$
\end{enumerate}

We have the following diagram.
\begin{equation*}
\begin{diagram}[small] H_{1}\left( E\left( K\left( \eta ,J\right) \right)
;\mathbb{Z}\Gamma _{n}\right) & \pile{\lTo^g\\ \rTo_{f_{\ast}}} &
H_{1}\left( E\left( K\right) ;\mathbb{Z}\Gamma _{n}\right) \\ \dTo>{B^{-1}}
& & \dTo>{B^{-1}} \\ H_{2}\left( E\left( K\left( \eta ,J\right) \right)
;\mathcal{K}_{n}/\mathbb{Z}\Gamma _{n}\right) & \rTo_{f_{\ast}} &
H_{2}\left( E\left( K\right) ;\mathcal{K}_{n}/\mathbb{Z}\Gamma _{n}\right)
\\ \dTo>{\pi} & & \dTo>{\pi} \\ H_{2}\left( E\left( K\left( \eta ,J\right)
\right) ,\partial E\left( K\left( \eta ,J\right) \right)
;\mathcal{K}_{n}/\mathbb{Z}\Gamma _{n}\right) & \rTo_{f_{\ast}} &
H_{2}\left( E\left( K\right) ,\partial E\left( K\right)
;\mathcal{K}_{n}/\mathbb{Z}\Gamma _{n}\right) \\ \dTo>{P.D.} & & \dTo>{P.D.}
\\ \overline{H^{1}\left( E\left( K\left( \eta ,J\right) \right)
;\mathcal{K}_{n}/\mathbb{Z}\Gamma _{n}\right) } & \lTo_{f^{\ast}} &
\overline{H^{1}\left( E\left( K\right) ;\mathcal{K}_{n}/\mathbb{Z}\Gamma
_{n}\right)} \\ \dTo>{\kappa} & & \dTo>{\kappa} \\ H_{1}\left( E\left(
K\left( \eta ,J\right) \right) ;\mathbb{Z}\Gamma _{n}\right) ^{\#} &
\pile{\rTo^{g^{\#}}\\ \lTo_{f^{\ast}}} & H_{1}\left( E\left( K\right)
;\mathbb{Z}\Gamma _{n}\right) ^{\#} \\ \end{diagram}
\end{equation*}
By the naturality of the Bockstein isomorphism, $\pi $, the
Poincar\'{e} duality isomorphism, and the Kronecker map, $f^{\ast
}\circ \mathcal{B\ell } _{K}\circ f_{\ast }=\mathcal{B\ell
}_{K\left( \eta ,J\right) }$.  So $ g^{\#}\circ f^{\ast }\circ
\mathcal{B\ell }_{K}\circ f_{\ast }\circ g=g^{\#}\circ
\mathcal{B\ell }_{K\left( \eta ,J\right) }\circ g$.  Since $ f_{\ast
}\circ g=id$ and $g^{\#}\circ f^{\ast }=id$, it follows that $
g^{\#}\circ \mathcal{B\ell }_{K\left( \eta ,J\right) }\circ
g=\mathcal{B\ell }_{K}$.  Hence the first claim is proved.

Consider the following diagram.
\begin{equation*}
\begin{diagram}[small] H_{1}\left( E\left( J\right) ;\mathbb{Z}\Gamma _{n}\right) &
\rTo^{i_{\ast}} & H_{1}\left( E\left( K\left( \eta ,J\right) \right)
;\mathbb{Z}\Gamma _{n}\right) \\ \dTo>{B^{-1}} & & \dTo>{B^{-1}} \\
H_{2}\left( E\left( J\right) ;\mathcal{K}_{n}/\mathbb{Z}\Gamma _{n}\right) &
\rTo^{i_{\ast}} & H_{2}\left( E\left( K\left( \eta ,J\right) \right)
;\mathcal{K}_{n}/\mathbb{Z}\Gamma _{n}\right) \\ \dTo>{\pi} & & \dTo>{\pi}
\\ H_{2}\left( E\left( J\right) ,\partial E\left( J\right)
;\mathcal{K}_{n}/\mathbb{Z}\Gamma _{n}\right) & & H_{2}\left( E\left(
K\left( \eta ,J\right) \right) ,\partial E\left( K\left( \eta ,J\right)
\right) ;\mathcal{K}_{n}/\mathbb{Z}\Gamma _{n}\right) \\ \dTo>{P.D.} & &
\dTo>{P.D.} \\ \overline{H^{1}\left( E\left( J\right)
;\mathcal{K}_{n}/\mathbb{Z}\Gamma _{n}\right) } & & \overline{H^{1}\left(
E\left( K\left( \eta ,J\right) \right) ;\mathcal{K}_{n}/\mathbb{Z}\Gamma
_{n}\right)} \\ \dTo>{\kappa} & & \dTo>{\kappa} \\ H_{1}\left( E\left(
J\right) ;\mathbb{Z}\Gamma _{n}\right) ^{\#} & \lTo^{i^{\ast}} & H_{1}\left(
E\left( K\left( \eta ,J\right) \right) ;\mathbb{Z}\Gamma _{n}\right) ^{\#}
\\ \end{diagram}
\end{equation*}
By the naturality of the Bockstein homomorphism, $B^{-1}\circ
i_{\ast }=i_{\ast }\circ B^{-1}.$ \ Consider the intersection
pairing (see, for example, \cite{duval}) $I_{E\left( K\left( \eta
,J\right) \right) }:H_{2}\left( E\left( K\left( \eta ,J\right)
\right) ;\mathcal{K}_{n}/ \mathbb{Z}\Gamma _{n}\right) \rightarrow
H_{1}\left( E\left( K\left( \eta ,J\right) \right) ;\mathbb{Z}\Gamma
_{n}\right) ^{\#}$ on $E\left( K\left( \eta ,J\right) \right) $
given by $I_{E\left( K\left( \eta ,J\right) \right) }=\kappa \circ
P.D.\circ \pi $.  Similarly, we have the intersection form
$I_{E\left( J\right) }=\kappa \circ P.D.\circ \pi :H_{2}\left(
E\left( J\right) ;\mathcal{K}_{n}/\mathbb{Z}\Gamma _{n}\right)
\rightarrow H_{1}\left( E\left( J\right) ;\mathbb{Z}\Gamma
_{n}\right) ^{\#}$ on $ E\left( J\right) $.  Since $i:E\left(
J\right) \rightarrow E\left( K\left( \eta ,J\right) \right) $ is an
embedding, $I_{E\left( J\right) }\left( x,y\right) =I_{E\left(
K\left( \eta ,J\right) \right) }\left( i_{\ast }\left( x\right)
,i_{\ast }\left( y\right) \right) $.  Therefore, $\kappa \circ
P.D.\circ \pi =i^{\ast }\circ \left( \kappa \circ P.D.\circ \pi
\right) \circ i_{\ast }$.  Thus,
\begin{eqnarray*}
\mathcal{B\ell }_{K\left( \eta ,J\right) }^{\otimes } &=&\left( \kappa \circ
P.D.\circ \pi \right) \circ B^{-1} \\
&=&i^{\ast }\circ \left( \kappa \circ P.D.\circ \pi \right) \circ i_{\ast
}\circ B^{-1} \\
&=&i^{\ast }\circ \left( \kappa \circ P.D.\circ \pi \right) \circ
B^{-1}\circ i_{\ast } \\
&=&i^{\ast }\circ \mathcal{B\ell }_{K\left( \eta ,J\right) }\circ i_{\ast }
\end{eqnarray*}
Therefore the second claim is proved.

Finally consider the following diagram.
\begin{equation*}
\begin{diagram} H_{1}\left( E\left( J\right) ;\mathbb{Z}\Gamma _{n}\right) &
\rTo^{i_{\ast }} & H_{1}\left( E\left( K\left( \eta ,J\right) \right)
;\mathbb{Z}\Gamma _{n}\right) & \pile{\lTo^g\\ \rTo_{f_{\ast}}} &
H_{1}\left( E\left( K\right) ;\mathbb{Z}\Gamma _{n}\right) \\
\dTo>{\mathcal{B\ell }_{K\left( \eta ,J\right) }^{\otimes }} & &
\dTo>{\mathcal{B\ell }_{K\left( \eta ,J\right) }} & & \dTo>{\mathcal{B\ell
}_{K}} \\ H_{1}\left( E\left( J\right) ;\mathbb{Z}\Gamma _{n}\right)^{\#} &
\lTo^{i^{\ast }} & H_{1}\left( E\left( K\left( \eta ,J\right) \right)
;\mathbb{Z}\Gamma _{n}\right)^{\#} & \pile{\rTo^{g^{\#}}\\ \lTo_{f^{\ast}}}
& H_{1}\left( E\left( K\right) ;\mathbb{Z}\Gamma _{n}\right)^{\#} \\
\end{diagram}
\end{equation*}
Since $f^{\ast }\circ \mathcal{B\ell }_{K}\circ f_{\ast
}=\mathcal{B\ell } _{K\left( \eta ,J\right) }$, it follows that
$g^{\#}\circ \mathcal{B\ell } _{K\left( \eta ,J\right) }\circ
i_{\ast }=g^{\#}\circ f^{\ast }\circ \mathcal{B\ell }_{K}\circ
f_{\ast }\circ i_{\ast }$.  But by Corollary \ref{SplitSES},
$f_{\ast }\circ i_{\ast }=0$.  Therefore $g^{\#}\circ \mathcal{
B\ell }_{K\left( \eta ,J\right) }\circ i_{\ast }=0$.  And since
$f_{\ast }\circ i_{\ast }=0$, we also have that $i^{\ast }\circ
f^{\ast }=0$.  Therefore, $i^{\ast }\circ \mathcal{B\ell }_{K\left(
\eta ,J\right) }\circ g=i^{\ast }\circ f^{\ast }\circ \mathcal{B\ell
}_{K}\circ f_{\ast }\circ g=0$.
\end{proof}

Recall that if $\eta \in \pi _{1}\left( E\left( K\right) \right)
^{\left( n\right) }$ and $\eta \notin \pi _{1}\left( E\left(
K\right) \right) ^{\left( n+1\right) }$ we have a ring monomorphism
$\psi: \mathbb{Z}\left[ t,t^{-1} \right] \rightarrow
\mathbb{Z}\Gamma _{n}$.  Therefore we have $\mathbb{Z}\left[
t,t^{-1}\right] $-module homomorphisms: $$\mathbb{Q} \left( t\right)
\hookrightarrow \mathcal{K}_{n}, \mathbb{Q}\left( t\right)
/\mathbb{Z}\left[ t,t^{-1}\right] \overset{\overline{\psi
}}{\hookrightarrow }\mathcal{K}_{n}/\mathbb{Z}\Gamma _{n}, H_{\ast
}\left( E\left( J\right) ;\mathbb{Z}\left[ t,t^{-1}\right] \right)
\overset{\psi _{\ast }}{ \rightarrow }H_{\ast }\left( E\left(
J\right) ;\mathbb{Z}\Gamma _{n}\right). $$ We will state the
following theorem without proof since the proof is quite technical
and it will not be needed for our main result.  It shows that
$\mathcal{B\ell }_{K\left( \eta ,J\right) }^{\otimes }$ is
determined by the classical Blanchfield form for $J$.

\begin{theorem}
\label{Tensor}If $\eta \in \pi _{1}\left( E\left( K\right) \right)
^{\left( n\right) }$ and $\eta \notin \pi _{1}\left( E\left(
K\right) \right) ^{\left( n+1\right) }$,
\begin{equation*}
\mathcal{B\ell }_{K\left( \eta ,J\right) }^{\otimes }\left( \psi
_{\ast }\left( x_{1}\right) ,\psi _{\ast }\left( x_{2}\right)
\right) =\overline{ \psi }\left( \mathcal{B\ell }_{J}\left(
x_{1},x_{2}\right) \right)
\end{equation*}
where $\mathcal{B\ell }_{J}$ is the classical Blanchfield linking
form for $ J $.
\end{theorem}

\section{Reducing from $\mathcal{B\ell }_{K\left( \protect\eta ,J\right) }$
to $\mathcal{B\ell }_{K\left( \protect\eta ,J\right) }^{\otimes }$}

Recall that our strategy is to find knots $J_{1}$ and $J_{2}$ such
that $ \mathcal{A}_{0}\left( J_{1}\right) \cong
\mathcal{A}_{0}\left( J_{2}\right) $ but $\mathcal{B\ell }_{K\left(
\eta ,J_{1}\right) }\ncong \mathcal{B\ell }_{K\left( \eta
,J_{2}\right) }$.  From Theorems \ref{DirectSum} and \ref{Tensor},
we know that if the classical Blanchfield linking forms of $J_{1}$
and $J_{2}$ are isomorphic, then $\mathcal{B\ell }_{K\left( \eta
,J_{1}\right) }^{\otimes }\cong \mathcal{B\ell }_{K\left( \eta
,J_{2}\right) }^{\otimes }$ and, therefore, $\mathcal{B\ell
}_{K\left( \eta ,J_{1}\right) }\cong \mathcal{B\ell }_{K\left( \eta
,J_{2}\right) }$.  However, the converses of these implications may
not follow.  That is, it may not be sufficient to choose $J_{1}$ and
$J_{2}$ with non-isomorphic classical Blanchfield linking forms.
 In this section, we find conditions on $K$ that
ensure that $\mathcal{B\ell }_{K\left( \eta ,J_{1}\right) }\cong
\mathcal{ B\ell }_{K\left( \eta ,J_{2}\right) }$ if and only if
$\mathcal{B\ell } _{K\left( \eta ,J_{1}\right) }^{\otimes }\cong
\mathcal{B\ell }_{K\left( \eta ,J_{2}\right) }^{\otimes }$.

\begin{proposition}
\label{ZGn'Torsion}If $\eta \in \pi _{1}\left( E\left( K\right)
\right) ^{\left( n\right) }$, $\eta \notin \pi _{1}\left( E\left(
K\right) \right) ^{\left( n+1\right) }$, and $n\geq 1$, then
$H_{1}\left( E\left( J\right) ; \mathbb{Z}\Gamma _{n}\right) $ is a
right $\mathbb{Z}\Gamma _{n}^{\prime }$-torsion module, where
$\Gamma _{n}^{\prime }=\left[\Gamma_{n},\Gamma _{n}\right] $. That
is, for any $\alpha \in H_{1}\left( E\left( J\right)
;\mathbb{Z}\Gamma _{n}\right) $, there is a non-zero $\gamma
_{\alpha }^{\prime }\in \mathbb{Z}\Gamma _{n}^{\prime }$ such that
$\alpha \gamma _{\alpha }^{\prime }=0$.
\end{proposition}

\begin{proof}
Recall from the proof of Corollary \ref{nckt8.2TimVersion} that $H_{1}\left(
E\left( J\right) ;\mathbb{Z}\Gamma _{n}\right) \cong \mathcal{A}_{0}\left(
J\right) \otimes _{\mathbb{Z}\left[ t,t^{-1}\right] }\mathbb{Z}\Gamma _{n}$.
 Hence, it suffices to consider $\beta \otimes \gamma $ where $\beta \in
\mathcal{A}_{0}\left( J\right) $ and $\gamma \in \mathbb{Z}\Gamma
_{n}$ are nonzero.  Let $\Delta _{J}$ be the classical Alexander
polynomial of $J$.  Since $\eta \in \pi _{1}\left( E\left( K\right)
\right) ^{\left( n\right) }\subset \pi _{1}\left( E\left( K\right)
\right) ^{\prime }$ and $\eta \notin \pi _{1}\left( E\left( K\right)
\right) ^{\left( n+1\right) }$, it follows that $\psi \left( \Delta
_{J}\right) \in \mathbb{Z}\Gamma _{n}^{\prime }$ is not zero. Since
$\mathbb{Z}\Gamma _{n}^{\prime }$ is a right $\mathbb{Z}\Gamma
_{n}-\left\{ 0\right\} $ Ore set \cite[p.16]{cohn}, there exist
$\widehat{\gamma }\in \mathbb{Z}\Gamma _{n}$ and $\gamma ^{\prime
}\in \mathbb{Z}\Gamma _{n}^{\prime }$ such that $\gamma \cdot \gamma
^{\prime }=\psi \left( \Delta _{J}\right) \cdot \widehat{\gamma }$
and $\widehat{\gamma }\neq 0$, $\gamma \neq 0$.  Thus
\begin{equation*}
\left( \beta \otimes \gamma \right) \cdot \gamma ^{\prime }=\beta \otimes
\left( \gamma \cdot \gamma ^{\prime }\right) =\beta \otimes \left( \psi
\left( \Delta _{J}\right) \cdot \widehat{\gamma }\right) =\left( \beta \cdot
\Delta _{J}\right) \otimes \widehat{\gamma }
\end{equation*}
since $\mathbb{Z}\Gamma _{n}$ is a left $\mathbb{Z}\left[
t,t^{-1}\right] $-module via the ring monomorphism $\psi $. However,
since $\beta \in \mathcal{A}_{0}\left( J\right) $ and $\Delta _{J}$
annihilates the Alexander module, $\beta \cdot \Delta _{J}=0$.
\end{proof}

\begin{lemma}
\label{Coefficients}$\mathbb{Z}\Gamma _{n}\cong \mathbb{Z}\pi
_{1}\left( E\left( K\right) \right) \otimes _{\mathbb{Z}\pi
_{1}\left( E\left( K\right) \right) ^{\prime }}\mathbb{Z}\Gamma
_{n}^{\prime }$ as right $\mathbb{Z} \Gamma _{n}^{\prime }$-modules.
\end{lemma}

\begin{proof}
Define $\varphi :\mathbb{Z}\Gamma _{n}\rightarrow \mathbb{Z}\pi
_{1}\left( E\left( K\right) \right) \otimes _{\mathbb{Z}\pi
_{1}\left( E\left( K\right) \right) ^{\prime }}\mathbb{Z}\Gamma
_{n}^{\prime }$ by $$\sum\limits_{i}n_{i} \left[ g_{i}\right]
\mapsto \left( \sum\limits_{i}n_{i}g_{i}\right) \otimes 1$$ where
$\left[ g\right] $ represents the coset of $g\in \pi _{1}\left(
E\left( K\right) \right) $ in $\Gamma _{n}=\pi _{1}\left( E\left(
K\right) \right) /\pi _{1}\left( E\left( K\right) \right) ^{\left(
n+1\right) }$.  First we must show that this is well-defined.  If $
h_{i}\in \pi _{1}\left( E\left( K\right) \right) ^{\left( n+1\right)
}$ for all $i$, then:
\begin{equation*}
\varphi \left( \sum\limits_{i}n_{i}\left[ g_{i}h_{i}\right] \right)
=\left( \sum\limits_{i}n_{i}g_{i}h_{i}\right) \otimes
1=\sum\limits_{i}n_{i}\left( g_{i}\otimes h_{i}\right)
=\sum\limits_{i}n_{i}\left( g_{i}\otimes 1\right).
\end{equation*}

For any $g\in \Gamma _{n}$ and $g^{\prime }\in \Gamma _{n}^{\prime
}$, $ gg^{\prime }\otimes 1=g\otimes g^{\prime }$. Hence $\varphi
\left( \left[ gg^{\prime }\right] \right) =\varphi \left( \left[
g\right] \right) g^{\prime }$.  It is easy to see that $\varphi $
preserves addition.  Finally, we define an inverse $\psi $ by
$$\sum\limits_{i}n_{i}g_{i}\otimes \sum\limits_{j}m_{j}g_{j}^{\prime
}\longmapsto \sum\limits_{i,j}n_{i}m_{j} \left[ g_{i}g_{j}^{\prime
}\right].$$  Note that $\psi \circ \varphi =id$ and $\varphi \circ
\psi =id$, since $$\sum\limits_{i}n_{i}g_{i}\otimes
\sum\limits_{j}m_{j}g_{j}^{\prime
}=\sum\limits_{i,j}n_{i}m_{j}g_{i}g_{j}^{\prime }\otimes 1.$$
Therefore $ \varphi $ is a right $\mathbb{Z}\Gamma _{n}^{\prime
}$-module isomorphism.
\end{proof}

\begin{proposition}[\protect\cite{nckt}, Prop. 9.3]
\label{nckt9.3}If $K$ is a fibered knot, then $\mathcal{A}_{n}\left(
K\right) $ has no $\mathbb{Z}\Gamma _{n}^{\prime }$-torsion.  That
is, for any $\beta \in \mathcal{A}_{n}\left( K\right) $ and $\gamma
^{\prime }\in \mathbb{Z}\Gamma _{n}^{\prime }$, if $\beta
\gamma^{\prime} =0$, then $\beta =0$ or $ \gamma ^{\prime }=0$.
\end{proposition}

\begin{proof}
Let $E\left( K\right) _{\infty }$ be the infinite cyclic cover, and
$\widetilde{E\left( K\right)}$ be the universal cover of $E\left(
K\right) $. Since $K$ is a fibered knot, $E\left( K\right) _{\infty
}$ is homotopy equivalent to a wedge of circles, $X$. Since $X$ is a
1-complex, $ H_{1}\left( X;\mathbb{Z}\Gamma _{n}^{\prime }\right)
\subset C_{1}\left( X; \mathbb{Z}\Gamma _{n}^{\prime }\right)$ which
is a free right $\mathbb{Z} \Gamma _{n}^{\prime }$-module. Therefore
$H_{1}\left( E\left( K\right) _{\infty };\mathbb{Z}\Gamma
_{n}^{\prime }\right) $ has no $\mathbb{Z}\Gamma _{n}^{\prime
}$-torsion.  Furthermore, the following are isomorphic right $
\mathbb{Z}\Gamma _{n}^{\prime }$-modules.
\begin{eqnarray*}
H_{1}\left( E\left( K\right) _{\infty };\mathbb{Z}\Gamma
_{n}^{\prime }\right) &=&H_{1}\left( C_{\ast }\left(\widetilde{
E\left( K\right) }\right) \otimes _{\mathbb{Z}\pi _{1}\left( E\left(
K\right) \right) ^{\prime }}
\mathbb{Z}\Gamma _{n}^{\prime }\right) \\
&\cong &H_{1}\left( C_{\ast }\left( \widetilde{E\left( K\right)
}\right) \otimes _{\mathbb{Z}\pi _{1}\left( E\left( K\right) \right)
}\mathbb{Z}\pi _{1}\left( E\left( K\right) \right) \otimes
_{\mathbb{Z}\pi _{1}\left(
E\left( K\right) \right) ^{\prime }}\mathbb{Z}\Gamma _{n}^{\prime }\right) \\
&&\hspace{2.5in}\text{ \cite[p.109]{hilton-stammbach}} \\
&\cong &H_{1}\left( C_{\ast }\left( \widetilde{E\left( K\right)
}\right) \otimes _{\mathbb{Z}\pi _{1}\left( E\left( K\right) \right)
}\mathbb{Z}\Gamma
_{n}\right) \text{ by Lemma \ref{Coefficients}} \\
&=&H_{1}\left( E\left( K\right) ;\mathbb{Z}\Gamma _{n}\right)
\end{eqnarray*}
Therefore $\mathcal{A}_{n}\left( K\right) $ has no $\mathbb{Z}\Gamma
_{n}^{\prime }$-torsion.
\end{proof}

\begin{theorem}
\label{Fibered}Let $K$ be a fibered knot. Suppose $\eta \in \pi
_{1}\left( E\left( K\right) \right) ^{\left( n\right) }$, \linebreak
$\eta \notin \pi _{1}\left( E\left( K\right) \right) ^{\left(
n+1\right) }$, and $n\geq 1$. If $\mathcal{B\ell }_{K\left( \eta
,J_{1}\right) }$ and $\mathcal{B\ell } _{K\left( \eta ,J_{2}\right)
}$ are isomorphic, then $\mathcal{B\ell } _{K\left( \eta
,J_{1}\right) }^{\otimes }$ and $\mathcal{B\ell }_{K\left( \eta
,J_{2}\right) }^{\otimes }$ are isomorphic.
\end{theorem}

\begin{remark}
Before proving the theorem, we remark that for any fibered knot $K$,
that is not the unknot, and any $n\geq 1$, there exists an $\eta $
such that $ \eta \in \pi _{1}\left( E\left( K\right) \right)
^{\left( n\right) }$ and $\eta \notin \pi _{1}\left( E\left(
K\right) \right) ^{\left( n+1\right) }$. \ This is because $K$ being
fibered implies that $\pi _{1}\left( E\left( K\right) \right)
^{\left( 1\right) }$ is isomorphic to a free group $F_{k}$ on $k$
generators ($k>1$, since $K$ is not the unknot).  Therefore,
$$\frac{ \pi _{1}\left( E\left( K\right) \right) ^{\left( n\right)
}}{\pi _{1}\left( E\left( K\right) \right) ^{\left( n+1\right) }}
\cong \frac{ F_{k}^{\left( n-1\right) }}{F_{k}^{\left( n\right)
}},$$ which is well-known to be non-trivial.
\end{remark}

\begin{proof}
Suppose $\mathcal{B\ell }_{K\left( \eta ,J_{1}\right) }$ and
$\mathcal{B\ell }_{K\left( \eta ,J_{2}\right) }$ are isomorphic
forms.  That is, there is a right $\mathbb{Z}\Gamma _{n}$-module
isomorphism $\psi :\mathcal{A} _{n}\left( K\left( \eta ,J_{1}\right)
\right) \rightarrow \mathcal{A} _{n}\left( K\left( \eta
,J_{2}\right) \right) $ such that for any $x,y\in
\mathcal{A}_{n}\left( K\left( \eta ,J_{1}\right) \right) $,
\begin{equation*}
\mathcal{B\ell }_{K\left( \eta ,J_{1}\right) }\left( x,y\right)
=\mathcal{ B\ell }_{K\left( \eta ,J_{2}\right) }\left( \psi \left(
x\right) ,\psi \left( y\right) \right) \text{.}
\end{equation*}
Using Theorem \ref{nckt8.2}, we have a right $\mathbb{Z}\Gamma
_{n}$-module isomorphism $\psi :H_{1}\left( E\left( K\right)
;\mathbb{Z}\Gamma _{n}\right) \oplus H_{1}\left( E\left(
J_{1}\right) ;\mathbb{Z}\Gamma _{n}\right) \rightarrow H_{1}\left(
E\left( K\right) ;\mathbb{Z}\Gamma _{n}\right) \oplus H_{1}\left(
E\left( J_{2}\right) ;\mathbb{Z}\Gamma _{n}\right) $. \ Since
$H_{1}\left( E\left( K\right) ;\mathbb{Z}\Gamma _{n}\right) $ has no
$\mathbb{Z}\Gamma _{n}^{\prime }$-torsion by Proposition
\ref{nckt9.3} and $H_{1}\left( E\left( J_{i}\right) ;\mathbb{Z}
\Gamma _{n}\right) $ is a $\mathbb{Z}\Gamma _{n}^{\prime }$-torsion
module, for $i=1,2$, by Proposition \ref{ZGn'Torsion}, it follows
that $\psi $ restricted to $H_{1}\left( E\left( J_{1}\right)
;\mathbb{Z}\Gamma _{n}\right) $ is a right $\mathbb{Z}\Gamma
_{n}$-module isomorphism between $ H_{1}\left( E\left( J_{1}\right)
;\mathbb{Z}\Gamma _{n}\right) $ and $ H_{1}\left( E \left(
J_{2}\right) ;\mathbb{Z}\Gamma_{n}\right)$.  Hence if $
x_{1},y_{1}\in H_{1}\left( E\left( J_{1}\right) ;\mathbb{Z}\Gamma
_{n}\right) $, then $\psi \circ i_{1}\left( x_{1}\right)
=i_{2}\left( x_{2}\right) $ and $\psi \circ i_{1}\left( y_{1}\right)
=i_{2}\left( y_{2}\right) $, for some $x_{2},y_{2}\in H_{1}\left(
E\left( J_{2}\right) ; \mathbb{Z}\Gamma _{n}\right) $.  Finally, we
have the following.
\begin{eqnarray*}
\mathcal{B\ell }_{K\left( \eta ,J_{1}\right) }^{\otimes }\left(
x_{1},y_{1}\right) &=&\mathcal{B\ell }_{K\left( \eta ,J_{1}\right)
}\left( i_{1}\left( x_{1}\right) ,i_{1}\left( y_{1}\right) \right)
\text{ by Theorem \ref{DirectSum}} \\
&=&\mathcal{B\ell }_{K\left( \eta ,J_{2}\right) }\left( \psi \left(
i_{1}\left( x_{1}\right) \right) ,\psi \left( i_{1}\left( y_{1}\right)
\right) \right) \\
&=&\mathcal{B\ell }_{K\left( \eta ,J_{2}\right) }\left( i_{2}\left(
x_{2}\right) ,i_{2}\left( y_{2}\right) \right) \\
&=&\mathcal{B\ell }_{K\left( \eta ,J_{2}\right) }^{\otimes }\left(
x_{2},y_{2}\right) \text{.}
\end{eqnarray*}
Therefore $\mathcal{B\ell }_{K\left( \eta ,J_{1}\right) }^{\otimes
}$ and $ \mathcal{B\ell }_{K\left( \eta ,J_{2}\right) }^{\otimes }$
are isomorphic forms.
\end{proof}

\section{Relating $\mathcal{B\ell }_{K\left( \protect\eta ,J\right)
}^{\otimes }$ to the Equivariant Intersection Form $\protect\lambda _{J}$}

Our main result is that there exist knots with isomorphic $n$th
higher-order Alexander modules, but non-isomorphic $n$th
higher-order linking forms. The idea behind the construction of our
examples is to infect the same knot $K$ along the same element $\eta
\in \pi _{1}\left( E\left( K\right) \right) ^{\left( n\right) }$,
$\eta \notin \pi _{1}\left( E\left( K\right) \right) ^{\left(
n+1\right) }$ by different knots $J_{1}$ and $J_{2}$ such that $
\mathcal{A}_{0}\left( J_{1}\right) \cong \mathcal{A}_{0}\left(
J_{2}\right) $.  Corollary \ref{nckt8.2TimVersion} implies that the
results of these infections will have isomorphic $i$-th order
Alexander modules for $0\leq i\leq n$.  From Theorem \ref{Fibered},
we know that if we choose $K$ to be a fibered knot, it suffices to
find examples of knots $J_{1}$ and $J_{2}$ such that
$\mathcal{A}_{0}\left( J_{1}\right) \cong \mathcal{A}_{0}\left(
J_{2}\right) $, but $\mathcal{B\ell }_{K\left( \eta ,J_{1}\right)
}^{\otimes }\ncong \mathcal{B\ell }_{K\left( \eta ,J_{2}\right)
}^{\otimes }$.  In this section, we will relate to $\mathcal{B\ell
}_{K\left( \eta ,J\right) }^{\otimes }$ a new linking form
$\widehat{\mathcal{B\ell }}_{K\left( \eta ,J\right) }$ defined on
the 0-framed surgery on $S^{3}$ along $J$, which we will then relate
to the equivariant intersection form, $\lambda _{J}$, on a
particular 4-manifold, $W_{J}$, associated to $J$.

Let $M_{J}$ denote the closed 3-manifold resulting from 0-framed
surgery on $ S^{3}$ along $J$.  The kernel of $\pi _{1}\left(
E\left( J\right) \right) \rightarrow \pi _{1}\left( M_{J}\right) $
is normally generated by the longitude of $J$, which is in the
kernel of $\pi _{1}\left( E\left( J\right) \right) \rightarrow
\Gamma _{n}$.  Hence $\pi _{1}\left( E\left( J\right) \right)
\rightarrow \Gamma _{n}$ factors through $\pi _{1}\left(
M_{J}\right) $, inducing a $\Gamma _{n}$ coefficient system on
$M_{J}$.

\begin{proposition}
\label{MJHomology}If $\eta \in \pi _{1}\left( E\left( K\right)
\right) ^{\left( n\right) }$ and $\eta \notin \pi _{1}\left( E\left(
K\right) \right) ^{\left( n+1\right) }$, then $$H_{1}\left( E\left(
J\right) ;\mathbb{Z }\Gamma _{n}\right) \cong H_{1}\left(
M_{J};\mathbb{Z}\Gamma _{n}\right).$$
\end{proposition}

\begin{proof}
We consider the Mayer-Vietoris sequence with $\mathbb{Z}\Gamma _{n}$
coefficients for $M_{J}\cong E\left( J\right) \cup _{\partial
E\left( J\right) }\left( D^{2}\times S^{1}\right) $.
\begin{equation*}
H_{1}\left( M_{J}\right) \overset{\partial _{\ast }}{\rightarrow }
H_{0}\left( \partial E\left( J\right) \right) \overset{\left( \psi
_{1},\psi _{2}\right) }{\rightarrow }H_{0}\left( E\left( J\right)
\right) \oplus H_{0}\left( D^{2}\times S^{1}\right)
\end{equation*}
By Lemma \ref{nckt8.3}, $\psi _{1}:H_{0}\left( \partial E\left(
J\right) \right) \rightarrow H_{0}\left( E\left( J\right) \right) $
is an isomorphism.  Therefore $\partial _{\ast }:H_{1}\left(
M_{J}\right) \rightarrow H_{0}\left( \partial E\left( J\right)
\right) $ is the trivial map.  Thus we have the following.
\begin{equation*}
H_{1}\left( \partial E\left( J\right) \right) \overset{\left( \psi
_{1},\psi _{2}\right) }{\rightarrow }H_{1}\left( E\left( J\right)
\right) \oplus H_{1}\left( D^{2}\times S^{1}\right) \overset{j_{\ast
}+k_{\ast }}{ \rightarrow }H_{1}\left( M_{J}\right) \rightarrow 0
\end{equation*}
By Lemma \ref{nckt8.3}, we also have that $\psi _{1}:H_{0}\left(
\partial E\left( J\right) \right) \rightarrow H_{0}\left( E\left(
J\right) \right) $ is the trivial map.  Since $\eta \notin \pi
_{1}\left( E\left( K\right) \right) ^{\left( n+1\right) }$ the
generator of $\pi _{1}\left( D^{2}\times S^{1};\mathbb{Z}\right) $,
$\mu _{J}$, gets unwound in the $\Gamma _{n}$ -cover. Hence
$H_{1}\left( D^{2}\times S^{1}\right) =0$.  Therefore, $ j_{\ast
}:H_{1}\left( E\left( J\right) \right) \rightarrow H_{1}\left(
M_{J}\right) $ is an isomorphism.
\end{proof}

Notice that since $H_{1}\left( E\left( J\right) ;\mathbb{Z}\Gamma
_{n}\right) $ is a $\mathbb{Z}\Gamma _{n}$-torsion module by
Proposition \ref{ZGnTorsion}, it follows that $H_{1}\left(
M_{J};\mathbb{Z}\Gamma _{n}\right) $ is as well.  So the Bockstein
map, $B:H_{2}\left( M_{J}; \mathcal{K}_{n}/\mathbb{Z}\Gamma
_{n}\right) \rightarrow H_{1}\left( M_{J}; \mathbb{Z}\Gamma
_{n}\right) $, is an isomorphism by the same argument as in Section
\ref{SectionMarker}.

\begin{definition}
If $\eta \in \pi _{1}\left( E\left( K\right) \right) ^{\left(
n\right) }$ and $\eta \notin \pi _{1}\left( E\left( K\right) \right)
^{\left( n+1\right) }$, define $\widehat{\mathcal{B\ell }}_{K\left(
\eta ,J\right) }:H_{1}\left( M_{J};\mathbb{Z}\Gamma _{n}\right)
\rightarrow H_{1}\left( M_{J};\mathbb{Z} \Gamma _{n}\right) ^{\#}$,
to be the composition of the following maps:
\begin{equation*}
H_{1}\left( M_{J};\mathbb{Z}\Gamma _{n}\right)
\overset{B^{-1}}{\rightarrow } H_{2}\left(
M_{J};\mathcal{K}_{n}/\mathbb{Z}\Gamma _{n}\right) \overset{P.D.}
{\rightarrow }\overline{H^{1}\left(
M_{J};\mathcal{K}_{n}/\mathbb{Z}\Gamma _{n}\right) }\overset{\kappa
}{\rightarrow }H_{1}\left( M_{J};\mathbb{Z} \Gamma _{n}\right) ^{\#}
\end{equation*}
where $P.D.$ is the Poincar\'{e} duality isomorphism and $\kappa $ is the
Kronecker evaluation map.
\end{definition}

\begin{theorem}
\label{tensortohat}If $\eta \in \pi _{1}\left( E\left( K\right)
\right) ^{\left( n\right) }$ and $\eta \notin \pi _{1}\left( E\left(
K\right) \right) ^{\left( n+1\right) }$, then $\mathcal{B\ell
}_{K\left( \eta ,J\right) }^{\otimes }\cong \widehat{\mathcal{B\ell
}}_{K\left( \eta ,J\right) }$.  That is,
\begin{equation*}
\mathcal{B\ell }_{K\left( \eta ,J\right) }^{\otimes }\left(
x_{1},y_{1}\right) =\widehat{\mathcal{B\ell }}_{K\left( \eta
,J\right) }\left( j_{\ast }\left( x_{1}\right) ,j_{\ast }\left(
y_{1}\right) \right).
\end{equation*}
\end{theorem}

\begin{proof}
Since $j$ is an inclusion map, the proof is the same as that of the
second claim in the proof of Theorem \ref{DirectSum}. Note however,
that since $\partial M_{J}=\emptyset$, the map $\pi$ which appears
in that proof is unnecessary.
\end{proof}

Therefore, we have reduced our problem to finding examples of knots
$J_{1}$ and $J_{2}$ such that $\mathcal{A}_{0}\left( J_{1}\right)
\cong \mathcal{A} _{0}\left( J_{2}\right) $, but
$\widehat{\mathcal{B\ell }}_{K\left( \eta ,J_{1}\right) }\ncong
\widehat{\mathcal{B\ell }}_{K\left( \eta ,J_{2}\right) }$.  To
accomplish this, we will relate $\widehat{\mathcal{B\ell }}_{K\left(
\eta ,J\right) }$ to the equivariant intersection form, $\lambda
_{J}$, on a particular 4-manifold, $W_{J}$, associated to $J$.

Since the bordism group $\Omega _{3}\left( S^{1}\right) =0$, we can
choose a 4-manifold $W_{J}$ which bounds $M_{J}$ and such that $\pi
_{1}\left( W_{J}\right) \cong \mathbb{Z}$, generated by the meridian
of $J$. Furthermore, by adding copies of $\pm \mathbb{CP}^{2}$, we
can choose $W_{J}$ so that the signature of it is zero.  (See
\cite{cot2}.)  Since the kernel of $\pi _{1}\left( E\left( J\right)
\right) \rightarrow \pi _{1}\left( W_{J}\right) $ is $\pi _{1}\left(
E\left( J\right) \right) ^{\prime }$, which is in the kernel of $\pi
_{1}\left( E\left( J\right) \right) \rightarrow \Gamma _{n}$, it
follows that $\pi _{1}\left( E\left( J\right) \right) \rightarrow
\Gamma _{n}$ factors through $\pi _{1}\left( W_{J}\right) $,
defining an induced $\Gamma _{n}$ coefficient system on $ W_{J}$.
Let $\lambda _{J}:H_{2}\left( W_{J};\mathbb{Z}\Gamma _{n}\right)
\rightarrow \overline{Hom_{\mathbb{Z}\Gamma _{n}}\left( H_{2}\left(
W_{J}; \mathbb{Z}\Gamma _{n}\right) ,\mathbb{Z}\Gamma _{n}\right) }$
be the equivariant intersection form on $W_{J}$ with
$\mathbb{Z}\Gamma _{n}$ coefficients.  That is, $\lambda_J$ is the
composition of the following maps.
\begin{eqnarray*}
&&H_{2}\left( W_{J};\mathbb{Z}\Gamma _{n}\right) \overset{\pi
}{\rightarrow } H_{2}\left( W_{J},M_{J};\mathbb{Z}\Gamma _{n}\right)
\overset{P.D.}{
\rightarrow }\overline{H^{2}\left( W_{J};\mathbb{Z}\Gamma _{n}\right) } \\
&&\hspace{1in}\overset{\kappa }{\rightarrow }\overline{Hom_{\mathbb{Z}\Gamma
_{n}}\left( H_{2}\left( W_{J};\mathbb{Z}\Gamma _{n}\right) ,\mathbb{Z}\Gamma
_{n}\right) }
\end{eqnarray*}

We recall the following definitions from \cite{ranicki}.

\begin{definition}[\protect\cite{ranicki}, pp. 60-61,145,181,242]
Let $S$ be a right denominator set for a ring with involution, $R$,
and let $ \mathcal{M}$ be an $R$-module.  A symmetric form over $R$,
$\alpha : \mathcal{M}\rightarrow \overline{\Hom_{R}\left(
\mathcal{M},R\right) }$, is $S$-non-singular if $\alpha \otimes
id:\mathcal{M}\otimes _{R}RS^{-1}\rightarrow
\overline{\Hom_{R}\left( \mathcal{M} ,R\right) }\otimes _{R}RS^{-1}$
is an $RS^{-1}$-module isomorphism.
\end{definition}

\begin{definition}[\protect\cite{ranicki}, p. 243]
The boundary of an $S$-non-singular symmetric form over $R$, $\alpha
:\mathcal{M}\rightarrow \overline{\Hom_{R}\left( \mathcal{M}
,R\right) }\equiv \mathcal{M}^{\ast }$, is the non-singular (even)
symmetric linking form over $\left( R,S\right) $ defined by
\begin{eqnarray*}
\partial \alpha :\coker\alpha &\rightarrow &\overline{\Hom
\nolimits_{R}\left( \mathcal{M}^{\ast },RS^{-1}/R\right) } \\
x &\longmapsto &\left( y\longmapsto x\left( z\right) \cdot
s^{-1}\right)
\end{eqnarray*}
for any $x,y\in \overline{\Hom_{R}\left( \mathcal{M},R\right) }
\equiv \mathcal{M}^{\ast }$, where $z\in \mathcal{M}$, $s\in S$ are
chosen so that $ys=\alpha \left( z\right) $.
\end{definition}

If $\eta \in \pi _{1}\left( E\left( K\right) \right) ^{\left(
n\right) }$ and $\eta \notin \pi _{1}\left( E\left( K\right) \right)
^{\left( n+1\right) }$, we will show in Proposition
\ref{S-Nonsingular} that the equivariant intersection form $\lambda
_{J}$, as above, is a $\left( \mathbb{Z}\Gamma _{n}-\left\{
0\right\} \right) $-non-singular symmetric form over
$\mathbb{Z}\Gamma _{n}$.  Furthermore, we will show in Theorem
\ref{Boundary} that $\partial \lambda _{J}\cong \widehat{
\mathcal{B\ell }}_{K\left( \eta ,J\right) }$.  In Section
\ref{section7}, this will allow us to fit $\lambda _{J}$ and
$\widehat{ \mathcal{B\ell }}_{K\left( \eta ,J\right) }$ into an
exact sequence of Witt groups. The result will be that we can
distinguish $\widehat{\mathcal{B\ell }}_{K\left( \eta ,J_{1}\right)
}$ from $ \widehat{\mathcal{B\ell }}_{K\left( \eta ,J_{2}\right) }$
by using an invariant of $\lambda_{J_1}$ and $\lambda_{J_2}$.

\begin{proposition}
\label{S-Nonsingular}If $\eta \in \pi _{1}\left( E\left( K\right) \right)
^{\left( n\right) }$ and $\eta \notin \pi _{1}\left( E\left( K\right)
\right) ^{\left( n+1\right) }$, the equivariant intersection form $\lambda
_{J}$, as above, is a $\left( \mathbb{Z}\Gamma _{n}-\left\{ 0\right\}
\right) $-non-singular symmetric form over $\mathbb{Z}\Gamma _{n}$.
\end{proposition}

\begin{proof}
In order to prove the proposition, we require Lemmas \ref{TensorP}
and \ref{Kronecker}.

\begin{lemma}
\label{TensorP}$H_{p}\left( M_{J};\mathbb{Z}\Gamma _{n}\right) \cong
H_{p}\left( M_{J};\mathbb{Z}\left[ t,t^{-1}\right] \right) \otimes
_{\mathbb{ Z}\left[ t,t^{-1}\right] }\mathbb{Z}\Gamma _{n}$ and
$H_{p}\left( W_{J}; \mathbb{Z}\Gamma _{n}\right) \cong H_{p}\left(
W_{J};\mathbb{Z}\left[ t,t^{-1}\right] \right) \otimes
_{\mathbb{Z}\left[ t,t^{-1}\right] }\mathbb{Z }\Gamma _{n}$ as right
$\mathbb{Z}\Gamma _{n}$-modules.
\end{lemma}

\begin{proof}
If $\widetilde{M_{J}}$ is the universal cover of $M_{J}$, we have the
following.
\begin{eqnarray*}
H_{p}\left( M_{J};\mathbb{Z}\Gamma _{n}\right) &=&H_{p}\left( C_{\ast
}\left( \widetilde{M_{J}}\right) \otimes _{\mathbb{Z}\pi _{1}\left(
M_{J}\right) }\mathbb{Z}\Gamma _{n}\right) \\
&\cong &H_{p}\left( C_{\ast }\left( \widetilde{M_{J}}\right) \otimes
_{ \mathbb{Z}\pi _{1}\left( M_{J}\right) }\mathbb{Z}\left[
t,t^{-1}\right] \otimes _{\mathbb{Z}\left[ t,t^{-1}\right]
}\mathbb{Z}\Gamma _{n}\right)
\text{ \cite[p. 109]{hilton-stammbach}} \\
&\cong &H_{p}\left( C_{\ast }\left( \widetilde{M_{J}}\right) \otimes
_{ \mathbb{Z}\pi _{1}\left( M_{J}\right) }\mathbb{Z}\left[
t,t^{-1}\right] \right) \otimes _{\mathbb{Z}\left[ t,t^{-1}\right]
}\mathbb{Z}\Gamma _{n}\\
&&\hspace{0.1in}\text{since }\mathbb{Z}\Gamma _{n}\text{ is a free,
hence flat, }\mathbb{Z}[ t,t^{-1}] \text{ module \cite[Lem. 1.3]
{passman}} \\
&\cong &H_{p}\left( M_{J};\mathbb{Z}\left[ t,t^{-1}\right] \right)
\otimes _{ \mathbb{Z}\left[ t,t^{-1}\right] }\mathbb{Z}\Gamma _{n}
\end{eqnarray*}
A similar argument holds for $W_{J}$.
\end{proof}

\begin{lemma}
\label{Kronecker}$\kappa :H^{2}\left( W_{J};\mathbb{Z}\Gamma
_{n}\right) \rightarrow Hom_{\mathbb{Z}\Gamma _{n}}\left(
H_{2}\left( W_{J};\mathbb{Z} \Gamma _{n}\right) ,\mathbb{Z}\Gamma
_{n}\right) $ is a $\mathbb{Z}\Gamma _{n}$-module isomorphism.
\end{lemma}

\begin{proof}
Since $\pi _{1}\left( W_{J}\right) $ is generated by the meridian of
$J$, which is identified to the longitude of $\eta $ in $E\left(
K\left( \eta ,J\right) \right) $, and $\eta $ gets unwound in the
$\mathbb{Z}\Gamma _{n}$ -cover, $H_{1}\left( W_{J};\mathbb{Z}\Gamma
_{n}\right) =0$.  By analyzing the Universal Coefficient Spectral
Sequence \cite[Thm 2.3]{levine} we have the following exact
sequence.
\begin{eqnarray*}
Ext_{\mathbb{Z}\Gamma _{n}}^{2}\left( H_{0}\left(
W_{J};\mathbb{Z}\Gamma _{n}\right) ,\mathbb{Z}\Gamma _{n}\right)
\rightarrow H^{2}\left( W_{J}; \mathbb{Z}\Gamma _{n}\right)
\overset{\kappa }{\rightarrow }Hom_{\mathbb{Z} \Gamma _{n}}\left(
H_{2}\left( W_{J};\mathbb{Z}\Gamma _{n}\right) ,\mathbb{Z}
\Gamma _{n}\right) \\
 \rightarrow Ext_{\mathbb{Z}\Gamma _{n}}^{3}\left( H_{0}\left(
W_{J}; \mathbb{Z}\Gamma _{n}\right) ,\mathbb{Z}\Gamma _{n}\right)
\end{eqnarray*}
By Lemma \ref{TensorP}, $H_{0}\left( W_{J};\mathbb{Z}\Gamma
_{n}\right) \cong H_{0}\left( W_{J};\mathbb{Z}\left[ t,t^{-1}\right]
\right) \otimes _{ \mathbb{Z}\left[ t,t^{-1}\right]
}\mathbb{Z}\Gamma _{n}$.  Since the homological dimension of
$H_{0}\left( W_{J};\mathbb{Z}\left[ t,t^{-1}\right] \right) $ is
$1$, and since $\mathbb{Z}\Gamma _{n}$ is a free, and therefore
flat, $\mathbb{Z}\left[ t,t^{-1}\right] $-module, $H_{0}\left(
W_{J};\mathbb{ Z}\left[ t,t^{-1}\right] \right) \otimes
_{\mathbb{Z}\left[ t,t^{-1}\right] } \mathbb{Z}\Gamma _{n}$ also has
homological dimension $1$. \ Therefore $Ext_{ \mathbb{Z}\Gamma
_{n}}^{p}\left( H_{0}\left( W_{J};\mathbb{Z}\Gamma _{n}\right)
,\mathbb{Z}\Gamma _{n}\right) =0$ for $p=2,3$.  Hence $\kappa
:H^{2}\left( W_{J};\mathbb{Z}\Gamma _{n}\right) \rightarrow
Hom_{\mathbb{Z} \Gamma _{n}}\left( H_{2}\left(
W_{J};\mathbb{Z}\Gamma _{n}\right) ,\mathbb{Z} \Gamma _{n}\right) $
is a $\mathbb{Z}\Gamma _{n}$-module isomorphism.
\end{proof}

Having proven these lemmas, we continue our proof of Proposition
\ref{S-Nonsingular}. Since $\lambda _{J}\cong \kappa \circ P.D.\circ
\pi $, it remains to be shown that $\pi \otimes id$, $P.D.\otimes
id$, and $\kappa \otimes id$ are $ \mathcal{K}_{n}$-module
isomorphisms.

Since $\mathcal{K}_{n}$ is a flat $\mathbb{Z}\Gamma _{n}$-module,
the following is an exact sequence of $\mathcal{K}_{n}$-modules.
\begin{eqnarray*}
H_{2}\left( M_{J};\mathbb{Z}\Gamma _{n}\right) \otimes
_{\mathbb{Z}\Gamma _{n}}\mathcal{K}_{n}\rightarrow H_{2}\left(
W_{J};\mathbb{Z}\Gamma _{n}\right) \otimes _{\mathbb{Z}\Gamma
_{n}}\mathcal{K}_{n}\overset{\pi \otimes id}{\rightarrow
}H_{2}\left( W_{J},M_{J};\mathbb{Z}\Gamma
_{n}\right) \otimes _{\mathbb{Z}\Gamma _{n}}\mathcal{K}_{n} \\
\hspace{1in}\rightarrow H_{1}\left( M_{J};\mathbb{Z}\Gamma
_{n}\right) \otimes _{\mathbb{Z}\Gamma _{n}}\mathcal{K}_{n}
\end{eqnarray*}
By Propositions \ref{ZGnTorsion} and \ref{MJHomology}, $H_{1}\left(
M_{J};\mathbb{Z} \Gamma _{n}\right) \otimes _{\mathbb{Z}\Gamma
_{n}}\mathcal{K}_{n}=0$. Again by the flatness of $\mathcal{K}_{n}$,
$H_{2}\left( M_{J}; \mathbb{Z}\Gamma _{n}\right) \otimes
_{\mathbb{Z}\Gamma _{n}}\mathcal{K} _{n}\cong H_{2}\left(
M_{J};\mathcal{K}_{n}\right) $.  By Poincar\'{e} Duality and the
Universal Coefficient Theorem, $H_{2}\left( M_{J};\mathcal{K}
_{n}\right) \cong Hom_{\mathbb{Z}\Gamma _{n}}\left( H_{1}\left(
M_{J}; \mathbb{Z}\Gamma _{n}\right) ,\mathcal{K}_{n}\right) $.  But
since $ H_{1}\left( M_{J};\mathbb{Z}\Gamma _{n}\right) $ is a
torsion-module, $Hom_{ \mathbb{Z}\Gamma _{n}}\left( H_{1}\left(
M_{J};\mathbb{Z}\Gamma _{n}\right) , \mathcal{K}_{n}\right) =0$.
Hence $\pi \otimes id$ is a $\mathcal{K}_{n}$-module isomorphism.

By Poincar\'{e} Duality and Lemma \ref{Kronecker}, it follows that $
P.D.\otimes id$ and $\kappa \otimes id$ are $\mathcal{K}_{n}$-module
isomorphisms.  Therefore,
\begin{equation*}
\lambda _{J}\otimes id:H_{2}\left( W;\mathbb{Z}\Gamma _{n}\right)
\otimes _{ \mathbb{Z}\Gamma _{n}}\mathcal{K}_{n}\rightarrow
\overline{Hom_{\mathbb{Z} \Gamma _{n}}\left( H_{2}\left(
W;\mathbb{Z}\Gamma _{n}\right) ,\mathbb{Z} \Gamma _{n}\right)
\otimes _{\mathbb{Z}\Gamma _{n}}\mathcal{K}_{n}}
\end{equation*}
is a $\mathcal{K}_{n}$-module isomorphism.  Thus, $\lambda _{J}$ is
a $ \left( \mathbb{Z}\Gamma _{n}-\left\{ 0\right\} \right)
$-non-singular symmetric form over $\mathbb{Z}\Gamma _{n}$.
\end{proof}

\begin{theorem}
\label{Boundary}If $\eta \in \pi _{1}\left( E\left( K\right) \right)
^{\left( n\right) }$ and $\eta \notin \pi _{1}\left( E\left(
K\right) \right) ^{\left( n+1\right) }$, it follows that $\partial
\lambda _{J}\cong \widehat{ \mathcal{B\ell }}_{K\left( \eta
,J\right) }:H_{1}\left( M_{J};\mathbb{Z} \Gamma _{n}\right)
\rightarrow H_{1}\left( M_{J};\mathbb{Z}\Gamma _{n}\right) ^{\#}$.
\end{theorem}

\begin{proof}
By definition, $\partial \lambda _{J}$ is defined as follows
\begin{eqnarray*}
\partial \lambda _{J}:\coker\lambda _{J} &\rightarrow &\left(
\overline{\Hom\nolimits_{\mathbb{Z}\Gamma _{n}}\left( H_{2}\left(
W_{J};\mathbb{Z}\Gamma _{n}\right) ,\mathbb{Z}\Gamma _{n}\right)
}\right)
^{\#} \\
x &\longmapsto &\left( y\longmapsto x\left( z\right) \cdot \gamma
^{-1}\right)
\end{eqnarray*}
for any $x,y\in \overline{\Hom\nolimits_{\mathbb{Z}\Gamma
_{n}}\left( H_{2}\left( W_{J};\mathbb{Z}\Gamma _{n}\right)
,\mathbb{Z}\Gamma _{n}\right) }$, where $z\in H_{2}\left(
W_{J};\mathbb{Z}\Gamma _{n}\right) $ and $\gamma \in
\mathbb{Z}\Gamma _{n}-\left\{ 0 \right\}$ are chosen so that
$y\gamma =\lambda _{J}\left( z\right) $.

Consider the following commutative diagram.

\begin{equation*}
\begin{diagram}[small] H_{2}\left( W_{J};\mathbb{Z}\Gamma _{n}\right) & \rTo^{\pi}
& H_{2}\left( W_{J},M_{J};\mathbb{Z}\Gamma _{n}\right) \\ &
\rdTo(3,4)_{\lambda_{J}} & \dTo>{P.D.} \\ & & \overline{H^{2}\left(
W_{J};\mathbb{Z}\Gamma _{n}\right)} \\ & & \dTo>{\kappa} \\ & &
\overline{Hom_{\mathbb{Z}\Gamma _{n}}\left( H_{2}\left(
W_{J};\mathbb{Z}\Gamma _{n}\right) ,\mathbb{Z}\Gamma _{n}\right)} \\
\end{diagram}
\end{equation*}

Recall from Lemma \ref{Kronecker} that $P.D.$ and $\kappa $ are
$\mathbb{Z} \Gamma _{n}$-module isomorphisms.  Since $\lambda
_{J}=\kappa \circ P.D.\circ \pi $, it follows that $\coker\lambda
_{J}=\left( \kappa \circ P.D.\right) \left( \coker\pi \right) $.

We define $\psi $ to be the following.
\begin{eqnarray*}
\psi :\coker\pi &\rightarrow &H_{2}\left( W_{J},M_{J};\mathbb{Z}
\Gamma _{n}\right) ^{\#} \\
a &\longmapsto &\left( b\longmapsto \left[ \left( \kappa \circ P.D.\right)
\left( a\right) \right] \left( z\right) \cdot \gamma ^{-1}\right)
\end{eqnarray*}
for any $a,b\in H_{2}\left( W_{J},M_{J};\mathbb{Z}\Gamma _{n}\right)
$, where $z\in H_{2}\left( W_{J};\mathbb{Z}\Gamma _{n}\right) $,
$\gamma \in \mathbb{Z}\Gamma _{n}$ are chosen so that $\pi \left(
z\right) =b\gamma $.

\begin{lemma}
$\partial \lambda _{J}$ is isomorphic to $\psi $ under $\kappa \circ P.D.$.
 That is, $\psi =\left( \kappa \circ P.D.\right) ^{\#}\circ \partial
\lambda _{J}\circ \left( \kappa \circ P.D.\right) $.
\end{lemma}

\begin{proof}
In order to show this, we must show that the following diagram
commutes.
\begin{equation*}
\begin{diagram}[small] \coker\lambda _{J} & \lTo^{\kappa \circ P.D.} &
\coker\pi \\ \dTo<{\partial \lambda _{J}} & & \dTo>{\psi}
\\ \left( \overline{\Hom\nolimits_{\mathbb{Z}\Gamma
_{n}}\left(H_{2}\left( W_{J};\mathbb{Z}\Gamma _{n}\right)
,\mathbb{Z}\Gamma _{n}\right) }\right) ^{\#} & \rTo^{\left( \kappa
\circ P.D.\right) ^{\#}} & H_{2}\left( W_{J},M_{J};\mathbb{Z}\Gamma
_{n}\right) ^{\#} \\ \end{diagram}
\end{equation*}
Suppose $a,b\in H_{2}\left( W_{J},M_{J};\mathbb{Z}\Gamma _{n}\right)
$ and $ z\in H_{2}\left( W_{J};\mathbb{Z}\Gamma _{n}\right) $,
$\gamma \in \mathbb{Z} \Gamma _{n}$ are chosen so that $\pi \left(
z\right) =b\gamma $.  Let $ x=\left( \kappa \circ P.D.\right) \left(
a\right) $ and $y=\left( \kappa \circ P.D.\right) \left( b\right) $.
 Then $\lambda _{J}\left( z\right) =\left( \kappa \circ P.D.\circ
\pi \right) \left( z\right) =\left( \kappa \circ P.D.\right) \left(
b\gamma \right) =y\gamma $.  By definition, $\left[ \left( \partial
\lambda _{J}\right) \left( x\right) \right] \left( y\right) =x\left(
z\right) \gamma ^{-1}$.  Therefore,
\begin{eqnarray*}
\left[ \left[ \left( \kappa \circ P.D.\right) ^{\#}\circ \partial \lambda
_{J}\circ \left( \kappa \circ P.D.\right) \right] \left( a\right) \right]
\left( b\right) &=&\left[ \left( \partial \lambda _{J}\circ \kappa \circ
P.D.\right) \left( a\right) \right] \left( \left( \kappa \circ P.D.\right)
\left( b\right) \right) \\
&=&\left[ \partial \lambda _{J}\left( x\right) \right] \left( y\right) \\
&=&x\left( z\right) \gamma ^{-1} \\
&=&\left[ \left( \kappa \circ P.D.\right) \left( a\right) \right] \left(
z\right) \gamma ^{-1} \\
&=&\left[ \psi \left( a\right) \right] \left( b\right)
\end{eqnarray*}
Therefore $\partial \lambda _{J}$ is isomorphic to $\psi $.
\end{proof}

In order to relate $\psi $ to $\widehat{\mathcal{B\ell }}_{K\left(
\eta ,J\right) }$, it will be easier to work with a more algebraic
definition of $ \psi $.  We begin by proving the following lemma.

\begin{lemma}
$\pi ^{\#}:\overline{\Hom_{\mathbb{Z}\Gamma _{n}}\left( H_{2}\left(
W_{J},M_{J};\mathbb{Z}\Gamma _{n}\right) ,\mathcal{K}_{n}\right) }
\rightarrow \overline{\Hom_{\mathbb{Z}\Gamma _{n}}\left( H_{2}\left(
W_{J};\mathbb{Z}\Gamma _{n}\right) ,\mathcal{K}_{n}\right) }$ is an
isomorphism.
\end{lemma}

\begin{proof}
In the proof of Proposition \ref{S-Nonsingular}, we showed that $\pi
\otimes id:H_{2}\left( W_{J};\mathbb{Z}\Gamma _{n}\right) \otimes
_{\mathbb{Z}\Gamma _{n}}\mathcal{K}_{n}\rightarrow H_{2}\left(
W_{J},M_{J};\mathbb{Z}\Gamma _{n}\right) \otimes _{\mathbb{Z}\Gamma
_{n}}\mathcal{K}_{n}$ is a $\mathcal{K }_{n}$-module isomorphism. \
Since $\mathcal{K}_{n}$ is a flat $\mathbb{Z} \Gamma _{n}$-module,
by Poincar\'{e} Duality and the Universal Coefficient Theorem,
\begin{equation*}
\pi ^{\#}:\overline{\Hom\nolimits_{\mathbb{Z}\Gamma _{n}}\left(
H_{2}\left( W_{J},M_{J};\mathbb{Z}\Gamma _{n}\right)
,\mathcal{K}_{n}\right) }\rightarrow
\overline{Hom_{\mathbb{Z}\Gamma_{n}}\left( H_{2}\left(
W_{J};\mathbb{Z}\Gamma _{n}\right) ,\mathcal{K}_{n}\right) }
\end{equation*}
is an isomorphism.
\end{proof}

By definition, $\left[ \psi \left( a\right) \right] \left( b\right)
=\left[ \left( \kappa \circ P.D.\right) \left( a\right) \right]
\left( z\right) \cdot \gamma ^{-1}$ where $\pi \left( z\right)
=b\gamma $. \ So \linebreak $ \left[ \pi ^{\#}\left( \psi \left(
a\right) \right) \right] \left( z\right) = \left[ \psi \left(
a\right) \right] \left( \pi \left( z\right) \right) = \left[ \psi
\left( a\right) \right] \left( b\gamma \right) =\left[ \left( \kappa
\circ P.D.\right) \left( a\right) \right] \left( z\right) $.  Hence
we have that $\psi $ is the composition of the following maps.
\begin{eqnarray*}
&&H_{2}\left( W_{J},M_{J};\mathbb{Z}\Gamma _{n}\right)
\overset{P.D.}{ \rightarrow }\overline{H^{2}\left(
W_{J};\mathbb{Z}\Gamma _{n}\right) } \overset{\kappa }{\rightarrow
}\overline{\Hom\nolimits_{\mathbb{Z} \Gamma _{n}}\left( H_{2}\left(
W_{J};\mathbb{Z}\Gamma _{n}\right) ,\mathbb{Z}
\Gamma _{n}\right) } \\
&&\hspace{0.25in}\rightarrow \overline{\Hom\nolimits_{\mathbb{Z}
\Gamma _{n}}\left( H_{2}\left( W_{J};\mathbb{Z}\Gamma _{n}\right)
,\mathcal{K }_{n}\right) }\overset{\left( \pi ^{\#}\right)
^{-1}}{\rightarrow }\overline{ \Hom\nolimits_{\mathbb{Z}\Gamma
_{n}}\left( H_{2}\left( W_{J},M_{J};
\mathbb{Z}\Gamma _{n}\right) ,\mathcal{K}_{n}\right) } \\
&&\hspace{0.25in}\rightarrow H_{2}\left( W_{J},M_{J};\mathbb{Z}\Gamma
_{n}\right) ^{\#}
\end{eqnarray*}

Recall that since $\pi _{1}\left( W_{J}\right) $ is generated by the
meridian of $J$, which is identified to the longitude of $\eta $ in
$E\left( K\left( \eta ,J\right) \right) $, and $\eta $ gets unwound
in the $\mathbb{Z} \Gamma _{n}$-cover, it follows that $H_{1}\left(
W_{J};\mathbb{Z}\Gamma _{n}\right) =0$.  Therefore, the following is
an exact sequence.
\begin{equation*}
H_{2}\left( W_{J};\mathbb{Z}\Gamma _{n}\right) \overset{\pi
}{\rightarrow } H_{2}\left( W_{J},M_{J};\mathbb{Z}\Gamma _{n}\right)
\overset{\partial _{\ast }}{\rightarrow }H_{1}\left(
M_{J};\mathbb{Z}\Gamma _{n}\right) \rightarrow 0
\end{equation*}
Hence $\coker\pi \equiv \frac{H_{2}\left( W_{J},M_{J};\mathbb{Z}
\Gamma _{n}\right) }{\im\pi }=\frac{H_{2}\left(
W_{J},M_{J};\mathbb{Z} \Gamma _{n}\right) }{\ker \partial _{\ast
}}\cong H_{1}\left( M_{J};\mathbb{Z }\Gamma _{n}\right) $. \ The
following lemma will complete the proof of the theorem.

\begin{lemma}
$\psi $ is isomorphic to $\widehat{\mathcal{B\ell }}_{K\left( \eta
,J\right) }$. \ That is, $\psi =\partial _{\ast }^{\#}\circ
\widehat{\mathcal{B\ell }} _{K\left( \eta ,J\right) }\circ \partial
_{\ast }$.
\end{lemma}

\begin{proof}
Recall that $\widehat{\mathcal{B\ell }}_{K\left( \eta ,J\right)
}=B^{-1}\circ P.D.\circ \kappa $, where $B:H_{2}\left(
M_{J};\mathcal{K}_{n}/ \mathbb{Z}\Gamma _{n}\right) \rightarrow
$\linebreak $H_{1}\left( M_{J}; \mathbb{Z}\Gamma _{n}\right) $.
Alternatively, we have $\widehat{\mathcal{ B\ell }}_{K\left( \eta
,J\right) }=P.D.\circ C^{-1}\circ \kappa $ since the following is a
commutative diagram of $\mathbb{Z}\Gamma _{n}$-module isomorphisms.
\begin{equation*}
\begin{diagram} H_{2}\left( M_{J};\mathcal{K}_{n}/\mathbb{Z}\Gamma
_{n}\right) & \rTo^{B} & H_{1}\left( M_{J};\mathbb{Z}\Gamma _{n}\right) \\
\dTo>{P.D.} & & \dTo>{P.D.} \\ \overline{H^{1}\left(
M_{J};\mathcal{K}_{n}/\mathbb{Z}\Gamma _{n}\right) } & \rTo^{C} &
\overline{H^{2}\left( M_{J};\mathbb{Z}\Gamma _{n}\right) } \\ \end{diagram}
\end{equation*}

We must show that the following diagram commutes.
\begin{equation*}
\begin{diagram}[small] H_{2}\left( W_{J},M_{J};\mathbb{Z}\Gamma _{n}\right) &
\rTo^{\partial _{\ast }} & H_{1}\left( M_{J};\mathbb{Z}\Gamma _{n}\right) \\
\dTo>{P.D.} & & \dTo>{P.D.} \\ \overline{H^{2}\left( W_{J};\mathbb{Z}\Gamma
_{n}\right)} & \rTo^{} & \overline{H^{2}\left( M_{J};\mathbb{Z}\Gamma
_{n}\right)} \\ \dTo>{\kappa} \\ \overline{Hom_{\mathbb{Z}\Gamma _{n}}\left(
H_{2}\left( W_{J};\mathbb{Z}\Gamma _{n}\right) ,\mathbb{Z}\Gamma
_{n}\right)} & & \dTo>{C^{-1}} \\ \dTo>{} \\ \overline{Hom_{\mathbb{Z}\Gamma
_{n}}\left( H_{2}\left( W_{J};\mathbb{Z}\Gamma _{n}\right)
,\mathcal{K}_{n}\right)} & & \overline{H^{1}\left(
M_{J};\mathcal{K}_{n}/\mathbb{Z}\Gamma _{n}\right)} \\ \dTo>{\left( \pi
^{\#}\right) ^{-1}} \\ \overline{Hom_{\mathbb{Z}\Gamma _{n}}\left(
H_{2}\left( W_{J},M_{J};\mathbb{Z}\Gamma _{n}\right)
,\mathcal{K}_{n}\right)} & & \dTo>{\kappa} \\ \dTo>{} \\ H_{2}\left(
W_{J},M_{J};\mathbb{Z}\Gamma _{n}\right) ^{\#} & \lTo^{\left( \partial
_{\ast }\right) ^{\#}} & H_{1}\left( M_{J};\mathbb{Z}\Gamma _{n}\right)^{\#}
\\ \end{diagram}
\end{equation*}
The top box commutes by the naturality of Poincar\'{e} Duality.

The short exact sequence of $\mathbb{Z}\Gamma _{n}$ chain groups $
0\rightarrow C_{\ast }\left( M_{J}\right) \rightarrow C_{\ast
}\left( W_{J}\right) \rightarrow C_{\ast }\left( W_{J},M_{J}\right)
\rightarrow 0$ gives rise to the following exact sequences.

\begin{eqnarray*}
0 &\rightarrow &C_{\ast }\left( W_{J},M_{J}\right) ^{\ast }\rightarrow
C_{\ast }\left( W_{J}\right) ^{\ast }\rightarrow C_{\ast }\left(
M_{J}\right) ^{\ast } \\
0 &\rightarrow &C_{\ast }\left( W_{J},M_{J}\right) ^{\circ }\rightarrow
C_{\ast }\left( W_{J}\right) ^{\circ }\rightarrow C_{\ast }\left(
M_{J}\right) ^{\circ } \\
0 &\rightarrow &C_{\ast }\left( W_{J},M_{J}\right) ^{\#}\rightarrow C_{\ast
}\left( W_{J}\right) ^{\#}\rightarrow C_{\ast }\left( M_{J}\right) ^{\#}
\end{eqnarray*}
Here $\mathcal{M}^{\ast }\equiv \overline{\Hom\nolimits_{\mathbb{Z}
\Gamma _{n}}\left( \mathcal{M},\mathbb{Z}\Gamma _{n}\right) }$,
$\mathcal{M} ^{\circ }\equiv
\overline{\Hom\nolimits_{\mathbb{Z}\Gamma _{n}}\left(
\mathcal{M},\mathcal{K}_{n}\right) }$, and \linebreak $\mathcal{M
}^{\#}\equiv \overline{\Hom\nolimits_{\mathbb{Z}\Gamma _{n}}\left(
\mathcal{M},\mathcal{K}_{n}/\mathbb{Z}\Gamma _{n}\right) }$.

Also the short exact sequence $0\rightarrow \mathbb{Z}\Gamma
_{n}\rightarrow \mathcal{K}_{n}\rightarrow
\mathcal{K}_{n}/\mathbb{Z}\Gamma _{n}\rightarrow 0 $ gives rise to
the following exact sequences.
\begin{eqnarray*}
0 &\rightarrow &C_{\ast }\left( M_{J}\right) ^{\ast }\rightarrow C_{\ast
}\left( M_{J}\right) ^{\circ }\rightarrow C_{\ast }\left( M_{J}\right) ^{\#}
\\
0 &\rightarrow &C_{\ast }\left( W_{J}\right) ^{\ast }\rightarrow C_{\ast
}\left( W_{J}\right) ^{\circ }\rightarrow C_{\ast }\left( W_{J}\right) ^{\#}
\\
0 &\rightarrow &C_{\ast }\left( W_{J},M_{J}\right)
^{\ast}\rightarrow C_{\ast }\left( W_{J},M_{J}\right) ^{\circ
}\rightarrow C_{\ast }\left( W_{J},M_{J}\right) ^{\# }
\end{eqnarray*}
The lemma now follows from the commutativity of the following
diagram.
\begin{equation*}
\begin{diagram}[tight] C_{1}\left( W_{J},M_{J}\right) ^{\ast } & & \rTo & &
C_{1}\left( W_{J},M_{J}\right) ^{\circ } & & \rTo & & C_{1}\left(
W_{J},M_{J}\right) ^{\#} & & \\ & \rdTo & & & \vLine & \rdTo & & & \vLine &
\rdTo & \\ \dTo & & C_{2}\left( W_{J},M_{J}\right) ^{\ast } & \rTo & \HonV &
& C_{2}\left( W_{J},M_{J}\right) ^{\circ } & \rTo & \HonV & & C_{2}\left(
W_{J},M_{J}\right) ^{\#} \\ & & \dTo & & \dTo & & \dTo & & \dTo & & \\
C_{1}\left( W_{J}\right) ^{\ast } & \hLine & \VonH & \rTo & C_{1}\left(
W_{J}\right) ^{\circ } &\hLine & \VonH & \rTo & C_{1}\left( W_{J}\right)
^{\#} & & \dTo \\ & \rdTo & & & \vLine & \rdTo & & & \vLine & \rdTo & \\
\dTo & & C_{2}\left( W_{J}\right) ^{\ast } & \rTo & \HonV & & C_{2}\left(
W_{J}\right) ^{\circ } & \rTo & \HonV & & C_{2}\left( W_{J}\right) ^{\#} \\
& & \dTo & & \dTo & & \dTo & & \dTo & & \\ C_{1}\left( M_{J}\right) ^{\ast }
& \hLine & \VonH & \rTo & C_{1}\left( M_{J}\right) ^{\circ } &\hLine & \VonH
& \rTo & C_{1}\left( M_{J}\right) ^{\#} & & \dTo \\ & \rdTo & & & & \rdTo &
& & & \rdTo & \\ & & C_{2}\left( M_{J}\right) ^{\ast } & \rTo & & &
C_{2}\left( M_{J}\right) ^{\circ } & \rTo & & & C_{2}\left( M_{J}\right)
^{\#} \end{diagram}
\end{equation*}
\end{proof}
\end{proof}

\section{L-Theory and the L$^{\text{2}}$-Signature}\label{section7}

If $R$ is a ring with involution and $S$ is a right denominator set
for $R$, then from \cite[pp.172,274]{ranicki} we have the following
exact sequence of Witt groups.
\begin{equation*}
L\left( R\right) \rightarrow L_{S}\left( RS^{-1}\right)
\overset{\partial }{ \rightarrow }L\left( R,S\right)
\end{equation*}
Here $L\left( R\right) $ is the Witt group of nonsingular symmetric forms
over $R$; $L_{S}\left( RS^{-1}\right) $ is the Witt group of $S$-nonsingular
symmetric forms over $R$; and $L\left( R,S\right) $ is the Witt group of
non-singular symmetric linking forms over $\left( R,S\right) $.

Recall that we have reduced our problem to finding examples of knots
$J_{1}$ and $J_{2}$ such that $\mathcal{A}_{0}\left( J_{1}\right)
\cong \mathcal{A} _{0}\left( J_{2}\right) $, but
$\widehat{\mathcal{B\ell }}_{K\left( \eta ,J_{1}\right) }\ncong
\widehat{\mathcal{B\ell }}_{K\left( \eta ,J_{2}\right) }$.  If $\eta
\in \pi _{1}\left( E\left( K\right) \right) ^{\left( n\right) }$ and
$\eta \notin \pi _{1}\left( E\left( K\right) \right) ^{\left(
n+1\right) }$, we showed in Proposition \ref{S-Nonsingular}, that
$\lambda _{J}\in L_{\left( \mathbb{Z}\Gamma _{n}-\left\{ 0\right\}
\right) }\left( \mathcal{K}_{n}\right) $, and in Theorem
\ref{Boundary}, that $\partial \lambda _{J}\cong
\widehat{\mathcal{B\ell }}_{K\left( \eta ,J\right) }$.  Hence we
need an invariant defined on $L_{\left( \mathbb{Z}\Gamma
_{n}-\left\{ 0\right\} \right) }\left( \mathcal{K}_{n}\right) $ that
is trivial on the image of $L\left( \mathbb{Z}\Gamma _{n}\right) $.
 In this
section, we will find that the desired invariant is the reduced
L$^{\text{2} } $-signature.  Furthermore, we will observe that in
our case, the reduced L $^{\text{2}}$-signature of $\lambda _{J}$ is
dependent only on the Levine-Tristram signatures of $J$.  We refer
the reader to Section 5 of \cite{cot1} and Sections 2 and 5 of
\cite{cot2} for more details about the L$ ^{\text{2}}$-signature.

\begin{proposition}[\protect\cite{cot1}, Cor. 5.7, Prop. 5.12]
\label{L2signature}The L$^{\text{2}}$-signature $\sigma _{\Gamma
}^{\left( 2\right) }$ is a real valued homomorphism on the Witt
group of nonsingular symmetric forms over $\mathcal{K}_{n}$,
$L\left( \mathcal{K}_{n}\right) $. Furthermore, the
L$^{\text{2}}$-signature equals the ordinary signature $ \sigma
_{0}$ on the image of $L\left( \mathbb{Z}\Gamma _{n}\right) $.
\end{proposition}

Therefore $\sigma _{\Gamma }^{\left( 2\right) }-\sigma _{0}$ satisfies the
desired conditions for our invariant. \ So we have the following definition.

\begin{definition}
The reduced L$^{\text{2}}$-signature of $\left( W_{J},\Gamma
_{n}\right) $, is defined to be $$\sigma _{\Gamma }^{\left( 2\right)
}\left( \lambda \left( W_{J}\right) \right) -\sigma _{0}\left(
\lambda \left( W_{J}\right) \right) $$ where $\lambda \left(
W_{J}\right) :H_{2}\left( W_{J};\mathbb{Z}\Gamma _{n}\right)
\rightarrow \overline{Hom_{\mathbb{Z}\Gamma _{n}}\left( H_{2}\left(
W_{J};\mathbb{Z}\Gamma _{n}\right) ,\mathbb{Z}\Gamma _{n}\right) }$
is the equivariant intersection form on $W_{J}$ with
$\mathbb{Z}\Gamma _{n}$ coefficients.
\end{definition}

By Proposition \ref{L2signature}, the reduced L$^{\text{2}}$-signature is a
well-defined real valued homomorphism on $L\left( \mathcal{K}_{n}\right) $.
 Furthermore, the reduced L$^{\text{2}}$-signature is zero on the image of
image of $L\left( \mathbb{Z}\Gamma _{n}\right) $.  Therefore, it
suffices to choose knots $J_{1}$ and $J_{2}$ such that
$\mathcal{A}_{0}\left( J_{1}\right) \cong \mathcal{A}_{0}\left(
J_{2}\right) $, but with reduced L$ ^{\text{2}}$-signatures of
$\left( W_{J_{1}},\Gamma _{n}\right) $ and $ \left( W_{J_{2}},\Gamma
_{n}\right) $ that are not equal.

\begin{proposition}[\protect\cite{cot1}, Prop. 5.13]
\label{cot5.13}If $\phi :\pi _{1}\left( W\right) \rightarrow \Gamma $
factors through a subgroup $\Sigma $, then $\sigma _{\Gamma }^{\left(
2\right) }\left( \lambda \left( W\right) \right) =\sigma _{\Sigma }^{\left(
2\right) }\left( \lambda \left( W\right) \right) $.
\end{proposition}

\begin{corollary}
\label{corto5.13}If $\eta \in \pi _{1}\left( E\left( K\right)
\right) ^{\left( n\right) }$ and $\eta \notin \pi _{1}\left( E\left(
K\right) \right) ^{\left( n+1\right) }$, the reduced \linebreak
L$^{\text{2}}$ -signatures of $\left( W_{J},\Gamma _{n}\right) $ and
$\left( W_{J},\mathbb{Z }\right) $ are equal.
\end{corollary}

\begin{proof}
Recall that our coefficient system on $W$ is defined by $\varphi $
in the following commutative diagram.
\begin{equation*}
\begin{diagram} \pi _{1}\left( E\left( J\right) \right) & \rTo^{i_{\ast }} &
\pi _{1}\left( E\left( K\left( \eta,J\right) \right) \right) & \rTo^{f_{\ast
}} & \pi _{1}\left(E\left( K\right) \right) & \rTo^{\phi } & \Gamma _{n} \\
\dTo<{} & & & & & \ruTo(7,2)_{\varphi} \\ \pi_{1}\left( W_{J} \right) \\
\end{diagram}
\end{equation*}
Recall also that $\pi_{1} \left( W_{J} \right) \cong \mathbb{Z}$
generated by a meridian of $J$, which is identified in $E \left( K
\left( \eta,J \right) \right)$ to $\eta$. Since we are assuming that
$\eta \in \pi _{1}\left( E\left( K\right) \right) ^{\left( n\right)
}$ and $\eta \notin \pi _{1}\left( E\left( K\right) \right) ^{\left(
n+1\right) }$, $\varphi $ is a monomorphism.  Therefore, $\varphi
:\pi _{1}\left( W_{J}\right) \rightarrow \Gamma _{n}$ factors
through $\mathbb{Z}$. \ It follows from Proposition \ref{cot5.13}
that $\sigma _{\Gamma _{n}}^{\left( 2\right) }\left( \lambda \left(
W_{J}\right) \right) =\sigma _{\mathbb{Z}}^{\left( 2\right) }\left(
\lambda \left( W_{J}\right) \right) $.  Hence the reduced L
$^{\text{2}}$-signatures of $\left( W_{J},\Gamma _{n}\right) $ and
$\left( W_{J},\mathbb{Z}\right) $ are equal.
\end{proof}

\begin{proposition}[\protect\cite{cot1}, Prop 5.1, Lemma 5.9(4)]
\label{cotII5.1}If $\eta \in \pi _{1}\left( E\left( K\right) \right)
^{\left( n\right) }$ and $\eta \notin \pi _{1}\left( E\left(
K\right) \right) ^{\left( n+1\right) }$, the reduced
L$^{\text{2}}$-signature of $ \left( W_{J},\mathbb{Z}\right) $ is
equal to the integral of the Levine-Tristram signatures of $J$,
integrated over the circle of unit length.
\end{proposition}

Thus we have our main theorem:

\begin{theorem}
Given any $n\geq 1$, suppose $K$ is a fibered knot, that is not the
unknot, and choose $\eta $ such that $\eta \in \pi _{1}\left(
E\left( K\right) \right) ^{\left( n\right) }$ and $\eta \notin \pi
_{1}\left( E\left( K\right) \right) ^{\left( n+1\right) }$.  Let $J$
be a knot such that the integral of the Levine-Tristram signatures
of $J$, integrated over the circle of unit length, is non-zero. Then
$\mathcal{A}_{i}\left( K\left( \eta ,J\right) \right) \cong
\mathcal{A}_{i}\left( K\left( \eta ,-J\right) \right) $ for $0\leq
i\leq n$ and $\mathcal{B\ell }_{i}\left( K\left( \eta ,J\right)
\right) \cong \mathcal{B\ell }_{i}\left( K\left( \eta ,-J\right)
\right) $ for $0\leq i<n$, but $\mathcal{B\ell }_{n}\left( K\left(
\eta ,J\right) \right) \ncong \mathcal{B\ell }_{n}\left( K\left(
\eta ,-J\right) \right) $.
\end{theorem}

\begin{proof}
By Corollary \ref{nckt8.1b}, for $0\leq i\leq n-1$,
$\mathcal{A}_{i}\left( K\left( \eta ,J\right) \right) \allowbreak
\cong \mathcal{A}_{i}\left( K\right) \allowbreak \cong
\mathcal{A}_{i}\left( K\left( \eta ,-J\right) \right) $. \
Furthermore, by Corollary \ref{nckt8.2TimVersion}, since the
classical Alexander modules of $J$ and $-J$ are isomorphic,
$\mathcal{A} _{n}\left( K\left( \eta ,J\right) \right) \cong
\mathcal{A}_{n}\left( K\right) \oplus \left( \mathcal{A}_{0}\left(
J\right) \otimes _{\mathbb{Z} \left[ t,t^{-1}\right]
}\mathbb{Z}\Gamma _{n}\right) \cong \mathcal{A} _{n}\left( K\right)
\oplus \newline \left( \mathcal{A}_{0}\left( -J\right) \otimes
_{\mathbb{Z}\left[ t,t^{-1} \right] } \mathbb{Z}\Gamma _{n}\right)
\cong \mathcal{A}_{n}\left( K\left( \eta ,-J\right) \right) $.

By Theorem \ref{lowerlinkingforms}, for $0\leq i\leq n-1$,
$\mathcal{B\ell } _{i}\left( K\left( \eta ,J\right) \right) \cong
\mathcal{B\ell }_{i}\left( K\right) \cong \mathcal{B\ell }_{i}\left(
K\left( \eta ,-J\right) \right) $.  Now suppose $\mathcal{B\ell
}_{n}\left( K\left( \eta ,J\right) \right) $ and $\mathcal{B\ell
}_{n}\left( K\left( \eta ,-J\right) \right) $ are isomorphic. Since
$K$ is a non-trivial fibered knot, it follows from Theorem
\ref{Fibered}, that $\mathcal{B\ell }_{K\left( \eta ,J\right)
}^{\otimes }$ and $\mathcal{B\ell }_{K\left( \eta ,-J\right)
}^{\otimes }$ are isomorphic, and therefore $\widehat{\mathcal{B\ell
}}_{K\left( \eta ,J\right) }\cong \widehat{\mathcal{B\ell
}}_{K\left( \eta ,-J\right) }$ by Theorem \ref{tensortohat}. Hence,
Proposition \ref{L2signature} implies that the reduced
L$^{\text{2}}$-signature of $\left( W_{J},\Gamma _{n}\right) $ and
$\left( W_{-J},\Gamma _{n}\right) $ are equal.  It follows from
Proposition \ref{cotII5.1} that $\int\nolimits_{\omega \in
S^{1}}\sigma _{\omega }\left( J\right) d\omega
=\int\nolimits_{\omega \in S^{1}}\sigma _{\omega }\left( -J\right)
d\omega $.  However, for all $ \omega \in \mathbb{C}$, $\sigma
_{\omega }\left( -J\right) =-\sigma _{\omega }\left( J\right) $. By
assumption $\int\nolimits_{\omega \in S^{1}}\sigma _{\omega }\left(
J\right) d\omega \neq 0$.  Therefore we have reached a
contradiction.  Hence $\mathcal{B\ell }_{n}\left( K\left( \eta
,J\right) \right) \ncong \mathcal{B\ell }_{n}\left( K\left( \eta
,-J\right) \right) $.
\end{proof}

\section{Example}

Since the trefoil is a fibered knot (\cite{rolfsen}, p.327), we can
use the left-handed trefoil for $K$.  Furthermore, since the
integral of the Levine-Tristram signatures of the left-handed
trefoil, integrated over the circle of unit length, is
$\frac{4}{3}$, we can also use the left-handed trefoil for $J$.
Finally, since the trefoil is equivalent to its reverse, it follows
that $-J$ is the mirror-image of $J$.  That is, $-J$ is the
right-handed trefoil.  Finally, we must choose $\eta $ so that $\eta
\in \pi _{1}\left( E\left( K\right) \right) ^{\left( n\right) }$ but
$\eta \notin \pi _{1}\left( E\left( K\right) \right) ^{\left(
n+1\right) }$.  For the case when $n=1$, it suffices to choose $\eta
$ to be a curve which clasps a band of the standard Seifert surface.
 Therefore, we have the following construction.

\begin{figure}[H]
\centering
\includegraphics[scale=.57]{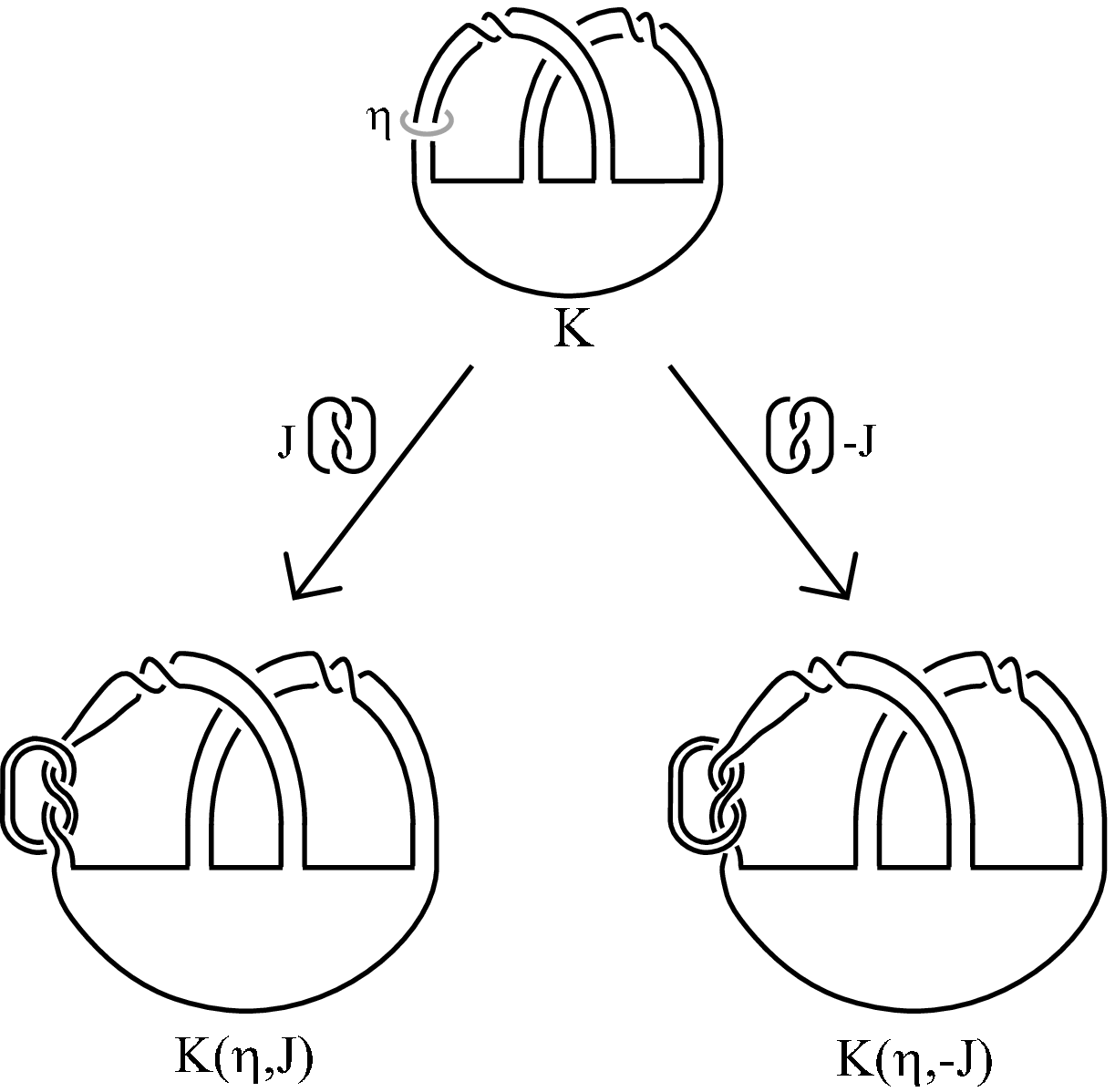}
\end{figure}


\begin{thebibliography}{COT1}
\bibitem[C]{nckt} T. Cochran, \textit{Noncommutative Knot Theory}, Alg.
Geom. Top., 4 (2004), 347-398.

\bibitem[COT1]{cot1} T. Cochran, K. Orr, and P. Teichner, \textit{Knot
Concordance, Whitney towers and L}$^{\text{2}}$\textit{-signature}, Annals
of Math., 157, (2003), 433-519.

\bibitem[COT2]{cot2} T. Cochran, K. Orr, and P. Teichner, \textit{Structure
in the Classical Knot Concordance Group}, Comment. Math. Helv. 79, 2004,
105-123.

\bibitem[Co]{cohn} P.M. Cohn, \textit{Skew Fields}, Cambridge University
Press, Cambridge, 1995.

\bibitem[DK]{daviskirk} J. Davis and P. Kirk, \textit{Lecture Notes in
Algebraic Topology}, American Mathematical Society, Providence, RI, 2001.

\bibitem[D]{duval} J. Duval, \textit{Forme de Blanchfield et cobordisme
d'entrelacs bords}, Comment. Math. Helv. 61, 1986, 617-635.

\bibitem[G]{camerongordon} C. Gordon, \textit{Some aspects of classical knot
theory}, Knot theory (Proc. Sem., Plans-sur-Bex, 1977), Springer, Berlin,
1978.

\bibitem[HS]{hilton-stammbach} P.J. Hilton and U. Stammbach, \textit{A
Course in Homological Algebra}, Springer-Verlag, New York, 1997.

\bibitem[Le]{levine} J.P. Levine, \textit{Knot Modules, I}, Trans. Amer.
Math. Soc. 229, 1977, 1-50.

\bibitem[Li]{lickorish} W.B.R. Lickorish, \textit{An Introduction to Knot
Theory}, Springer-Verlage, New York, 1997.

\bibitem[K]{kearton} C. Kearton, \textit{Blanchfield Duality and Simple Knots%
}, Trans. Amer. Math. Soc., 202 (1975), 141-160.

\bibitem[MR]{McConnellRobson} J.C. McConnell and J.C. Robson, \textit{%
Noncommutative Noetherian Rings}, Graduate Studies in Mathematics, 30.
American Mathematical Society, Providence, RI, 2001.

\bibitem[P]{passman} D. Passman, \textit{The Algebraic Structure of Group
Rings}, John Wiley and Sons, New York, 1977.

\bibitem[R]{ranicki} A. Ranicki, \textit{Exact Sequences in the Algebraic
Theory of Surgery}, Mathematical Notes, 26. Princeton University Press,
Princeton, NJ, 1981.

\bibitem[Ro]{rolfsen} D. Rolfsen, \textit{Knots and Links}, Publish or
Perish, Inc., Houston, TX, 1976.

\bibitem[Ste]{stenstrom} B. Stenstr\"{o}m, \textit{Rings of Quotients},
Springer-Verlag, New York, 1975.

\bibitem[T]{trotter} H. F. Trotter, \textit{On S-equivalence of Seifert
Matrices}, Invent. Math., 20, (1973), 173--207.
\end{thebibliography}
\end{document}